\documentclass[graybox]{svmult}

\usepackage{type1cm}        
\usepackage{makeidx}         
\usepackage{graphicx}        
\usepackage{multicol}        
\usepackage[bottom]{footmisc}

\usepackage{newtxtext}       %
\usepackage[varvw]{newtxmath}       

\makeindex             

\begin{document}

\title* {Natural damping of time-harmonic waves and its influence on Schwarz methods}
\author{Martin Gander\orcidID{0000-0001-8450-9223} and\\ Hui Zhang \orcidID{0000-0001-7245-0674}}
\institute{Martin Gander \at Department of Mathematics, University of Geneva, Rue du
  Général Dufour, 1211 Genève 4, Switzerland, \email{martin.gander@unige.ch} \and Hui Zhang
  (corresponding) \at Department of Applied Mathematics, Xi'an Jiaotong-Liverpool University, 111
  Ren'ai Road, 215123 Suzhou, China, \email{hui.zhang@xjtlu.edu.cn}}
%
%
\maketitle

\abstract*{The influence of various damping on the performance of
  Schwarz methods for time-harmonic waves is visualized by Fourier
  analysis.}

\section{Introduction}
\label{sec:1}

The Helmholtz equation $(\Delta+\omega^2)u=f$ is often used as the
prototype time harmonic problem which is difficult to solve by
numerical methods, both because high mesh resolutions are required
when the wave number $\omega$ becomes large and iterative solvers
struggle to solve such discretized problems \cite{ernst2011difficult}.
An idea is to use the shifted Helmholtz equation
\cite{van2007spectral}, $(\Delta+\omega^2+i\epsilon)u=f$, also called
the Helmholtz equation with damping, and the parameter $\epsilon$ can
be chosen to make iterative methods succeed, even though a difficult
compromise must be chosen between approximation of the Helmholtz
equation solution ($\epsilon\le O(\omega)$ \cite{gander2015applying})
and easy solvability by iterative methods ($\epsilon\ge O(\omega^2)$
\cite{cocquet2017large}) when using the shifted problem as a
preconditioner.

\begin{figure}
  \centering
  \mbox{\includegraphics[width=0.33\textwidth]{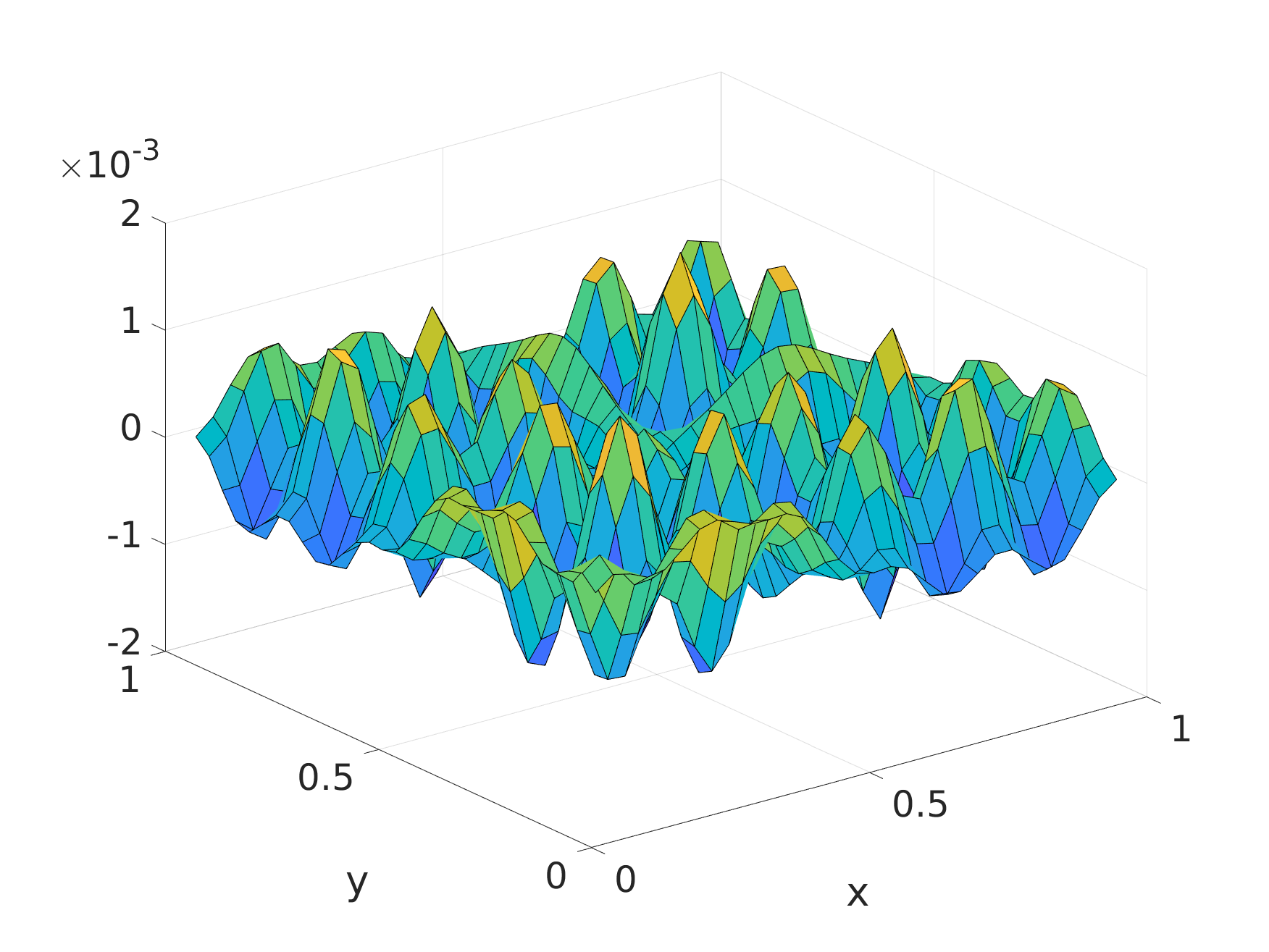}
  \includegraphics[width=0.33\textwidth]{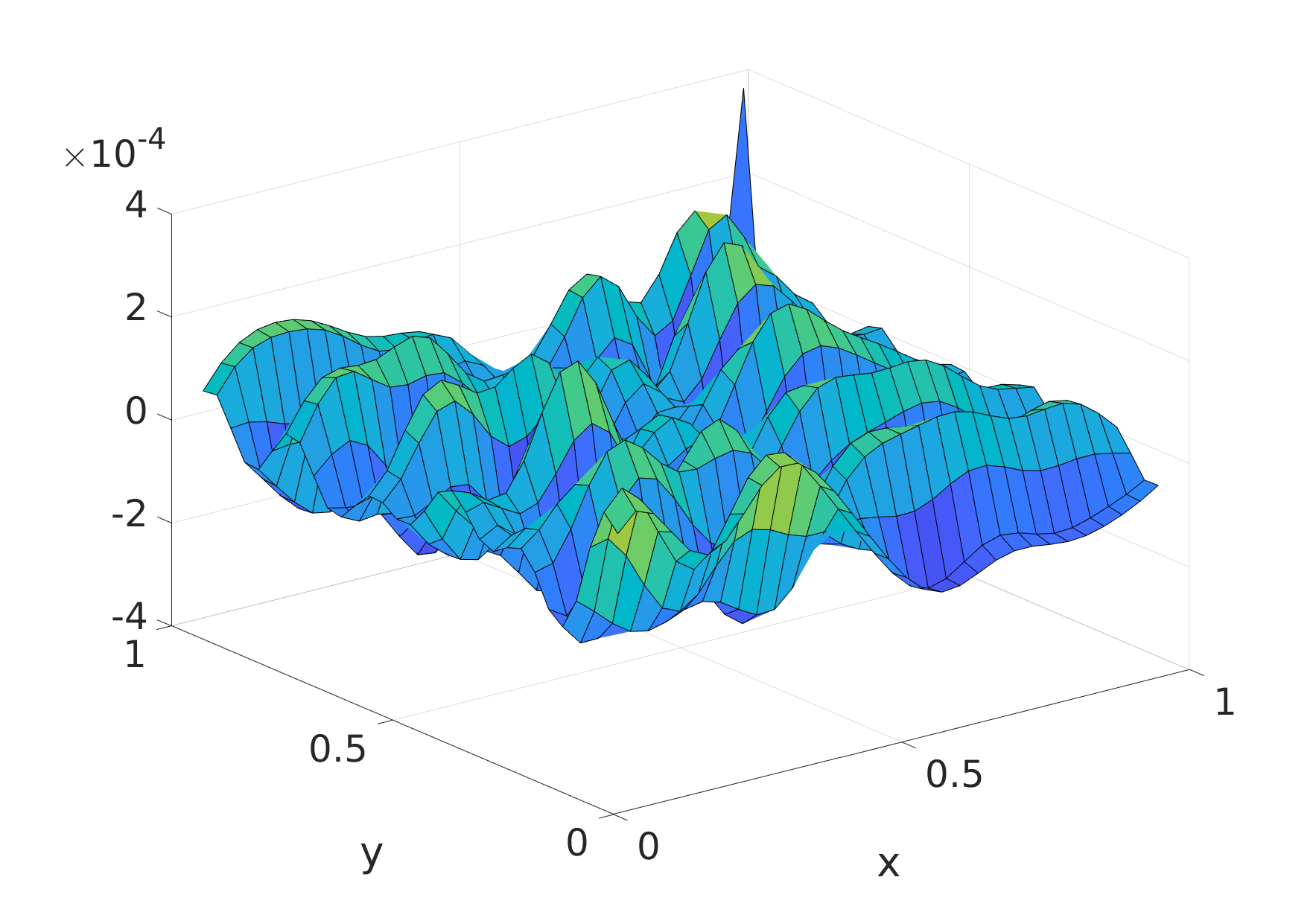}
  \includegraphics[width=0.33\textwidth]{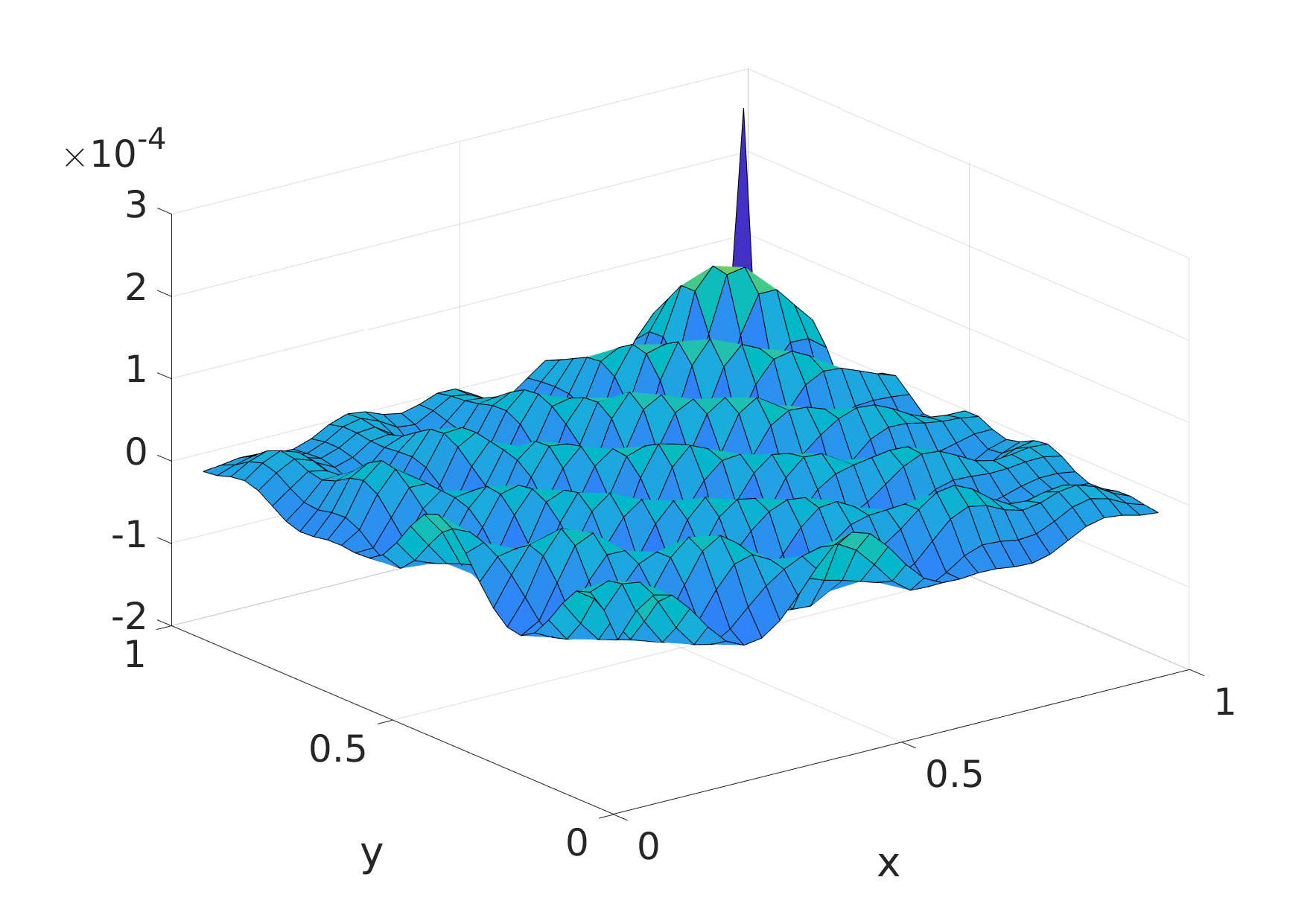}}
  \caption{Greens function for 3 Helmholtz problems. Left: closed
    cavity; middle: wave guide; right: free space, where optimized
    Schwarz solvers are as effective as for Laplace problems
    \cite{Gander_Zhang_2022}.}
  \label{FigSol}
\end{figure}

It was shown in \cite{Gander_Zhang_2022} that important damping is
also coming from the outer boundary conditions applied to the
Helmholtz equation, and while it is the closed cavity (all Dirichlet)
and the waveguide (Dirichlet along the guide and impedance at the
ends) are the really hard to solve by iterative methods, the free
space problem (e.g. impedance all around) becomes as easy to solve as the
Laplace problem in the constant coefficient case for Schwarz methods
of sweeping type \cite{gander2019class} based on optimized Schwarz
technology. We show in Fig. \ref{FigSol} a Greens function
corresponding to these three cases to illustrate how different the
solution looks like.

\begin{figure}
  \centering
  \mbox{\includegraphics[width=0.33\textwidth]{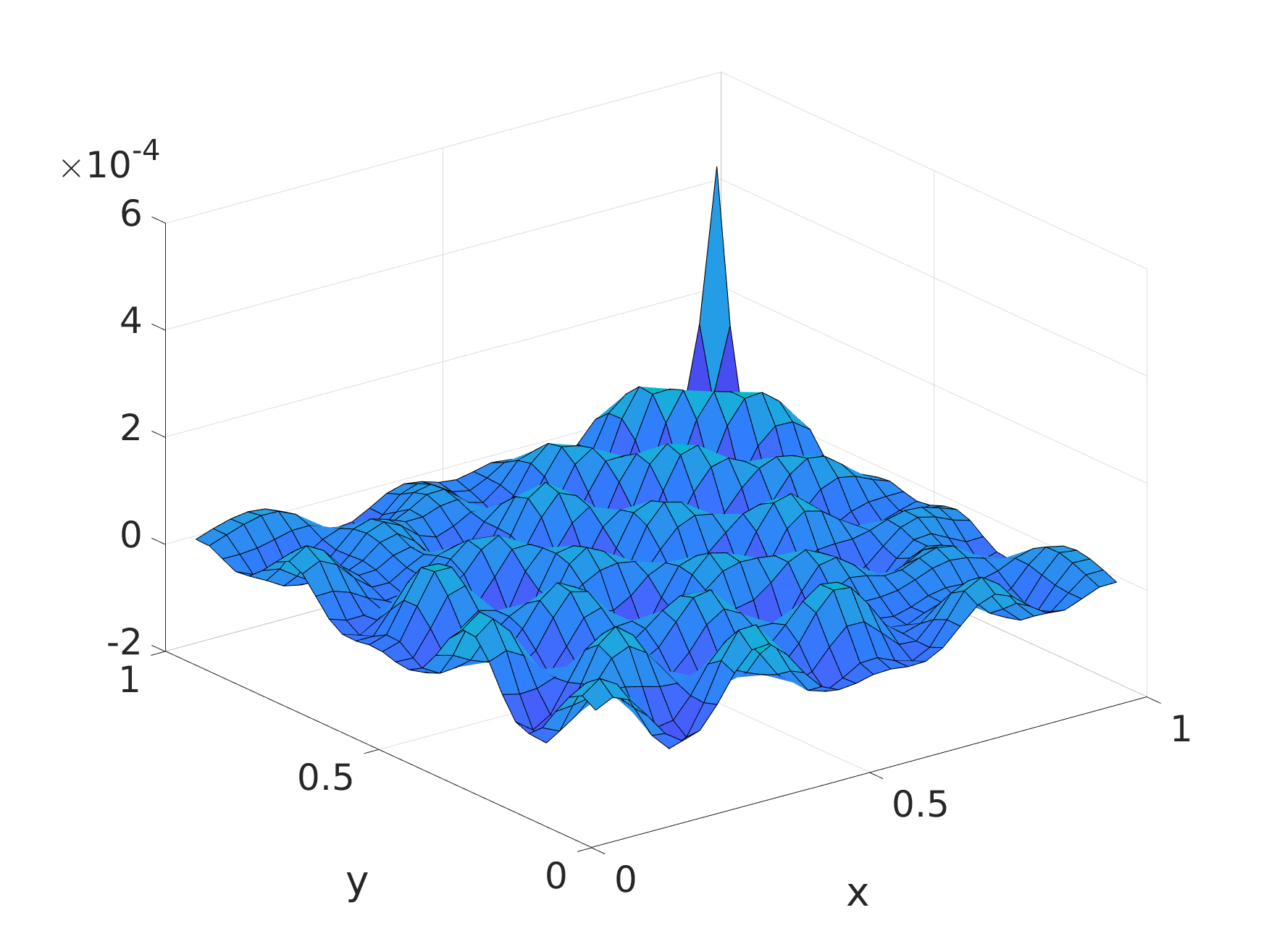}
  \includegraphics[width=0.33\textwidth]{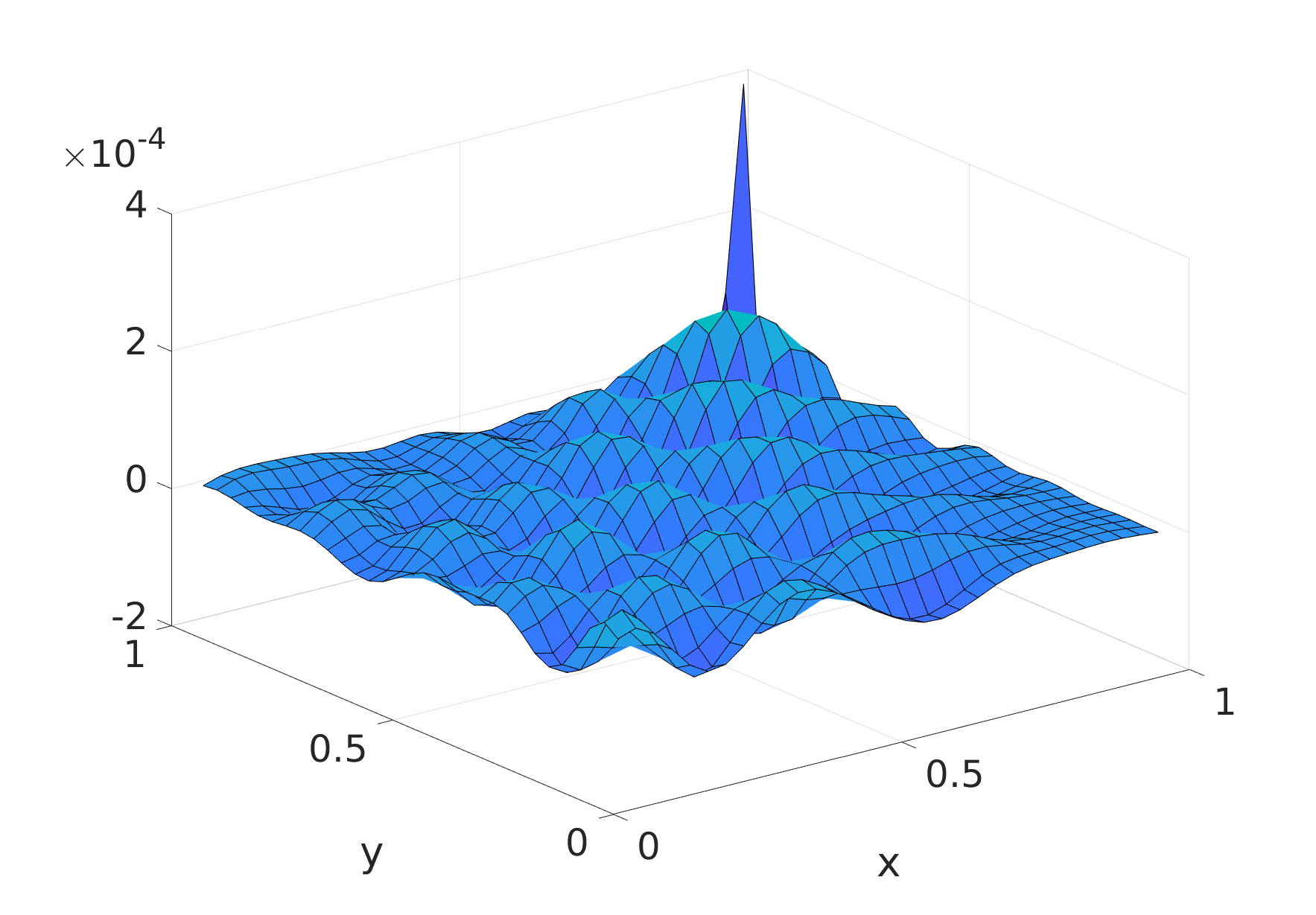}
  \includegraphics[width=0.33\textwidth]{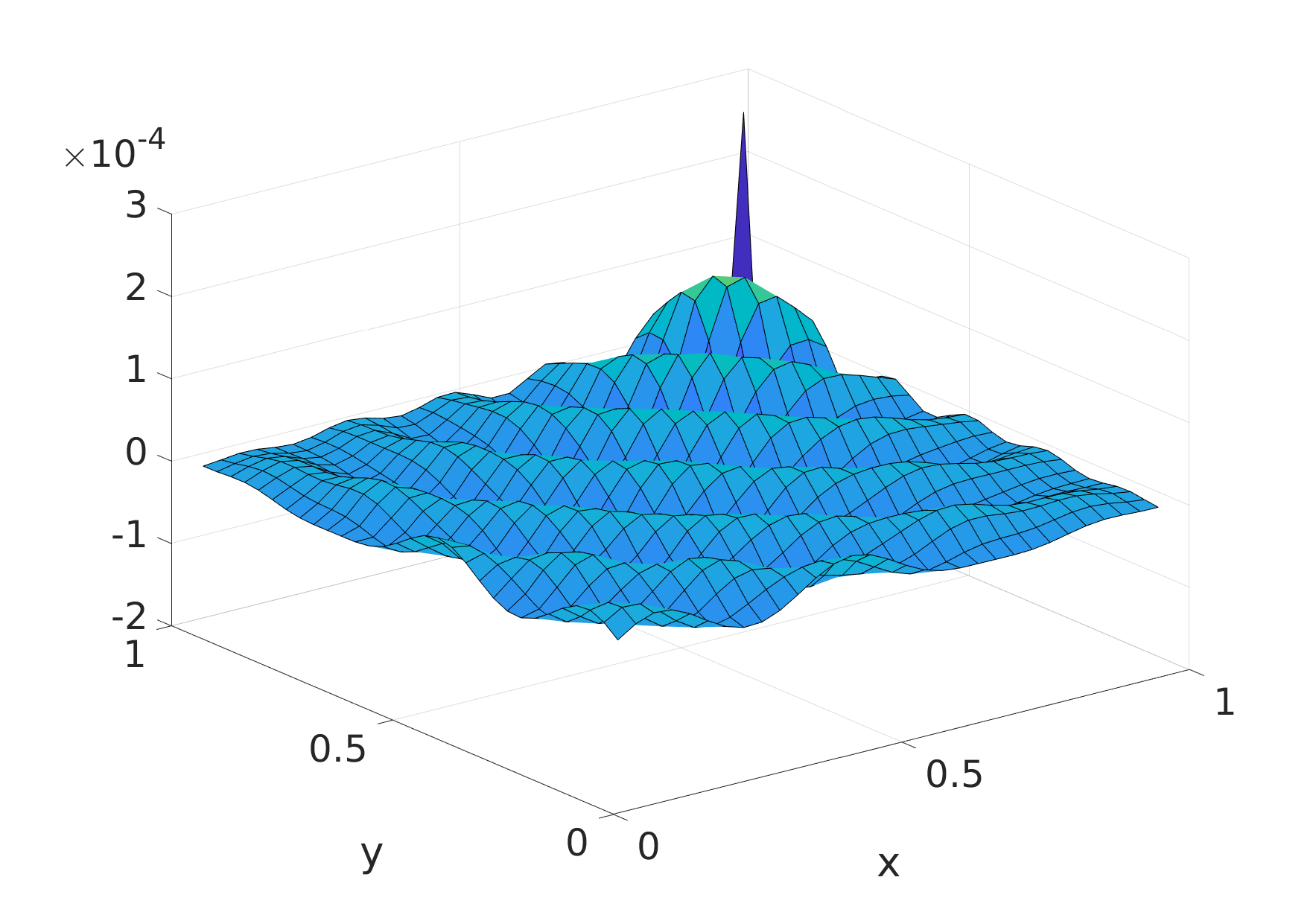}}
  \mbox{\includegraphics[width=0.33\textwidth]{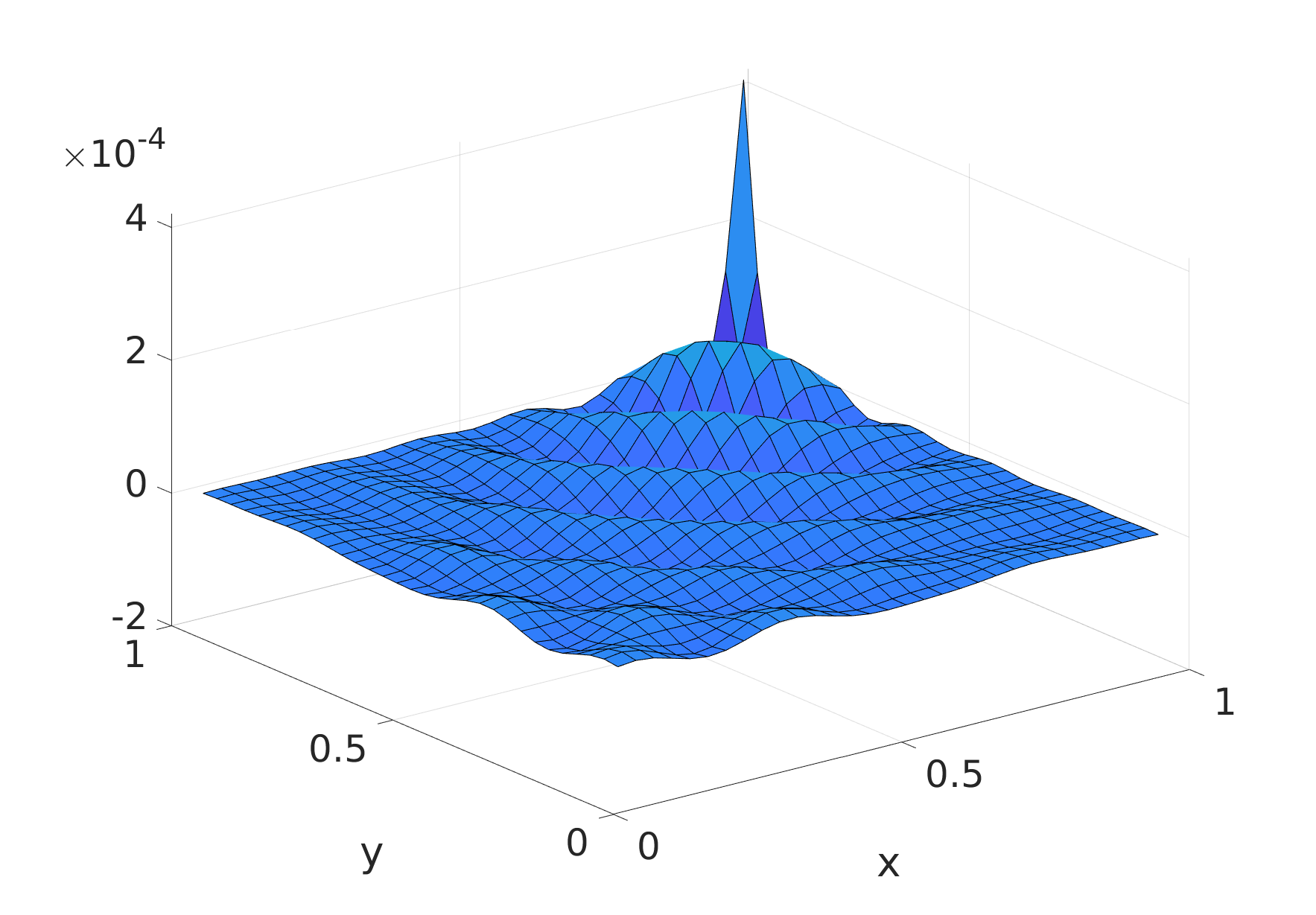}
  \includegraphics[width=0.33\textwidth]{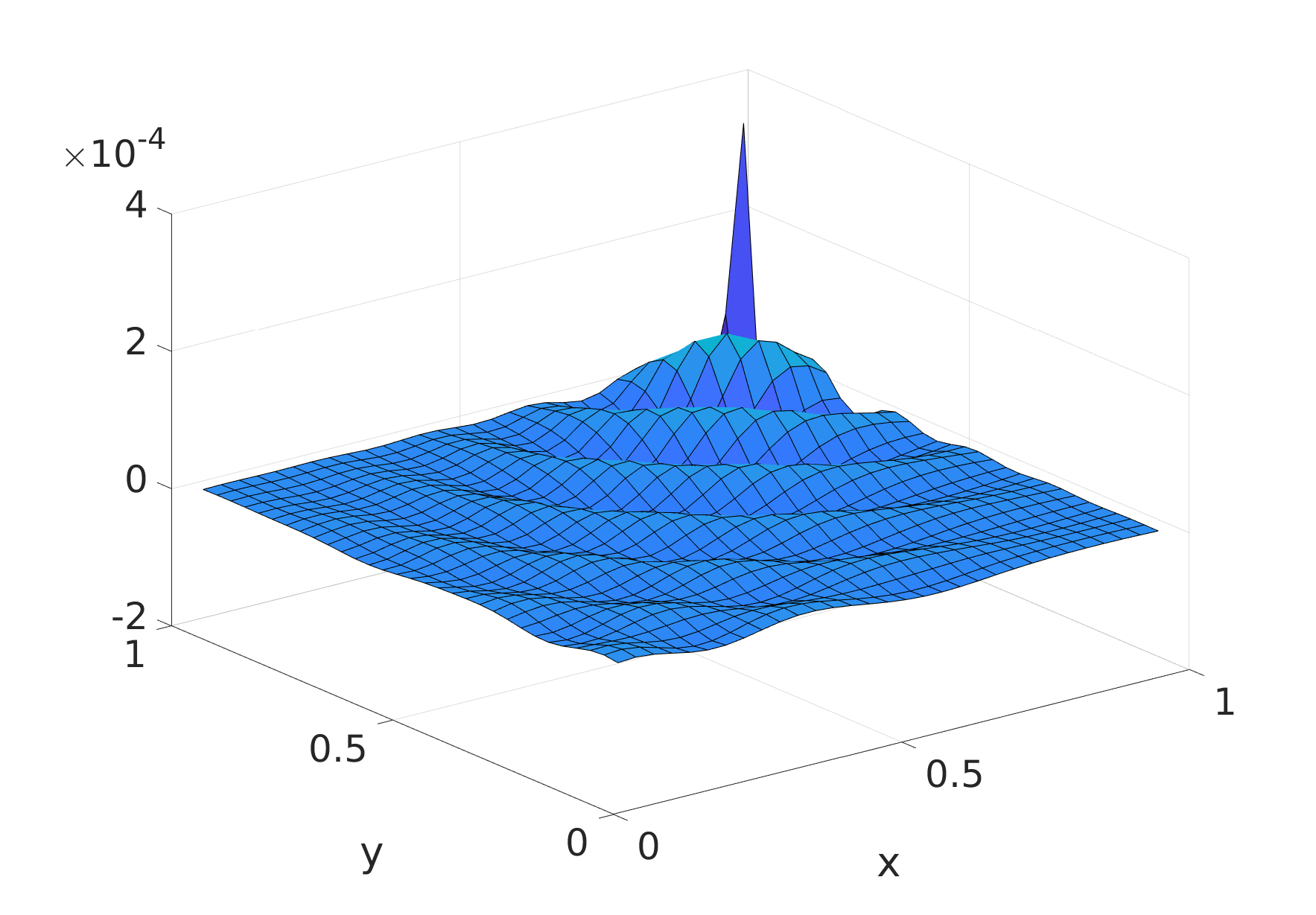}
  \includegraphics[width=0.33\textwidth]{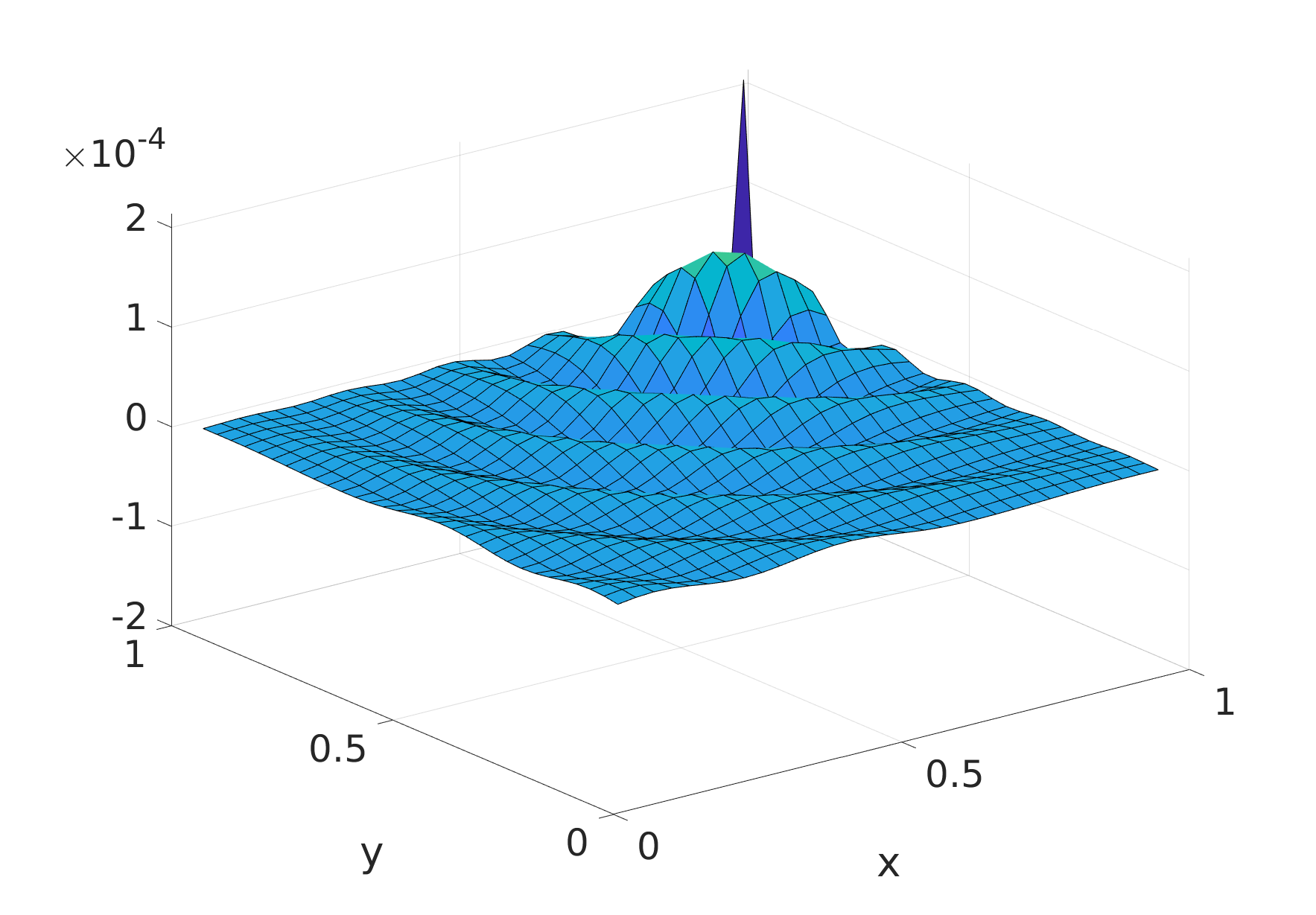}}
  \caption{Greens functions with first order damping $r=1$ (top) and
    viscoelastic damping $\gamma=0.003$ (bottom).}
  \label{FigSolr}
\end{figure}

We are interested here in studying if natural damping mechanisms
coming from the physical properties of the wave phenomenon one wants
to study can make the damped Helmholtz problem easy to solve by such
Schwarz methods. To do so, we go back to the underlying phyiscal wave
equation from which the Helmholtz equation arises, namely the second
order wave equation, $u_{tt}=\Delta u+f$.  There are two main damping
mechanisms that arise in nature for damping solutions of the second
order wave equation, first order and viscoelastic damping, 
\begin{equation}\label{WaveEquationWithPhysicalDamping}
  u_{tt}+ru_t=\Delta u+\gamma \partial_t\Delta u+f,
\end{equation}
where $r$ is the first order damping strength, and $\gamma$ is the
viscoelastic damping strength. These terms are present in all physical
phenomena in nature, the second order wave equation is a simplifying
model, and waves never remain forever, there is always damping.  If
the source $f(x,t)$ and the solution $u(x,t)$ have the time-harmonic
form $\hat{f}(x)\mathrm{e}^{\I\omega t}$ and
$\hat{u}(x)\mathrm{e}^{\I\omega t}$, then we find $(\omega^2 -
\I\omega r)\hat{u} + (1+\I\omega\gamma)\Delta \hat{u} = -\hat{f}$,
which we write it in the normalized form (omitting the hats) to obtain
the new Helmholtz equation that includes natural damping,
\begin{equation}\label{NewHelmholtz}
\Delta u- \eta u = -f/(1+\gamma\I\omega), \quad
  \eta:=-\omega^2(1-r\I\omega^{-1})(1+\gamma\I\omega)^{-1}.
\end{equation}
We show in Figure \ref{FigSolr} the influence of first order and
viscoelastic damping on the Greens functions from Figure \ref{FigSol}.
We see that both types of damping make the Greens function look like
the easy to solve case in Figure \ref{FigSol} on the right, all
difficulties from the closed cavity and wave guide configurations seem
to have disappeared. It is therefore of interest to investigate the
performance of Schwarz methods in these naturally damped
configurations.

We analyze here the parallel Schwarz method applied to
\eqref{NewHelmholtz} on the domain $\Omega:=(0,1)^2$ with impedance
transmission conditions at subdomain interfaces of the form
$\partial_{\mathbf{n}}u+\sqrt{\eta} u$, where $\mathbf{n}$ is the unit
normal. We decompose the domain along the $x$ direction into
overlapping subdomains of the same width $H+L$ with $L$ the overlap
width shared by two neighbors. We apply Fourier series in $y$, then
calculate the spectral radius of the iteration matrix of the interface
values $\partial_{\mathbf{n}}u+\sqrt{\eta} u$, which we call
convergence factor $\rho$, and present results both for the waveguide
and the closed cavity problems.

\section{Waveguide problem}
\label{sec2}

In the waveguide problem, $\partial_{\mathbf{n}}u+\sqrt{\eta} u=0$ at
$x\in\{0,1\}$, and $u=0$ on the other sides, and we consider first
order damping ($r>0$, $\gamma=0$) and viscoelastic damping ($r=0$,
$\gamma>0$) separately. We visualize the convergence factor $\rho$ in
the Fourier domain with $\xi$ the Fourier frequency for $y$.

\subsection{Helmholtz operator $\Delta + \omega^2 - \I\omega r$ in the waveguide}
\label{sec2.1}

\begin{figure}[t]
  \centering
  \mbox{\includegraphics[width=0.33\textwidth,trim=10 5 42 37,clip]{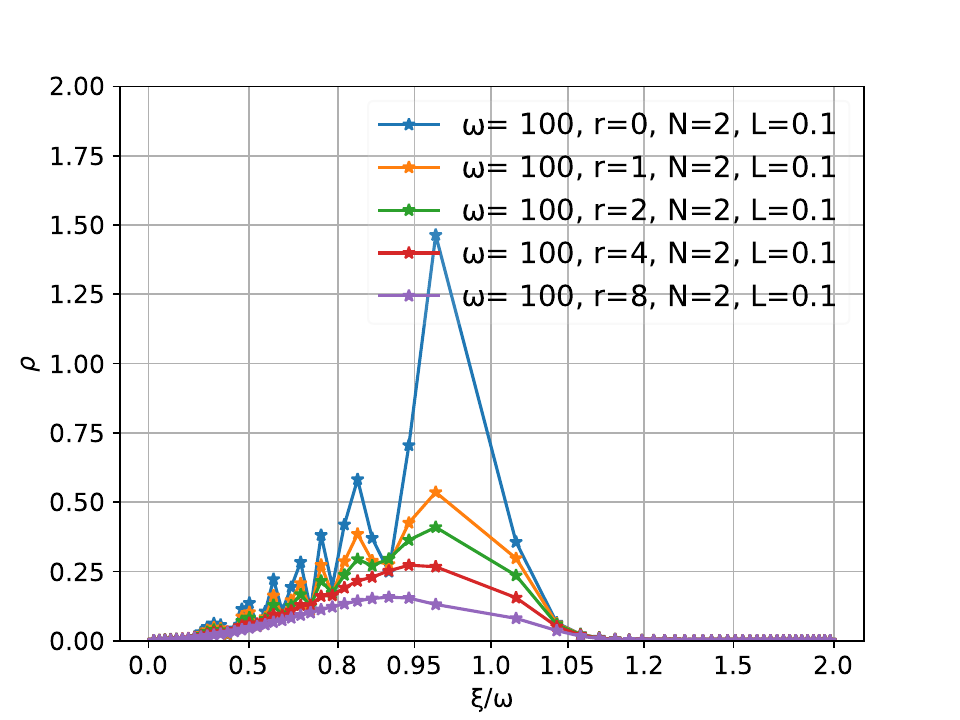}%
  \includegraphics[width=0.33\textwidth,trim=10 5 42 37,clip]{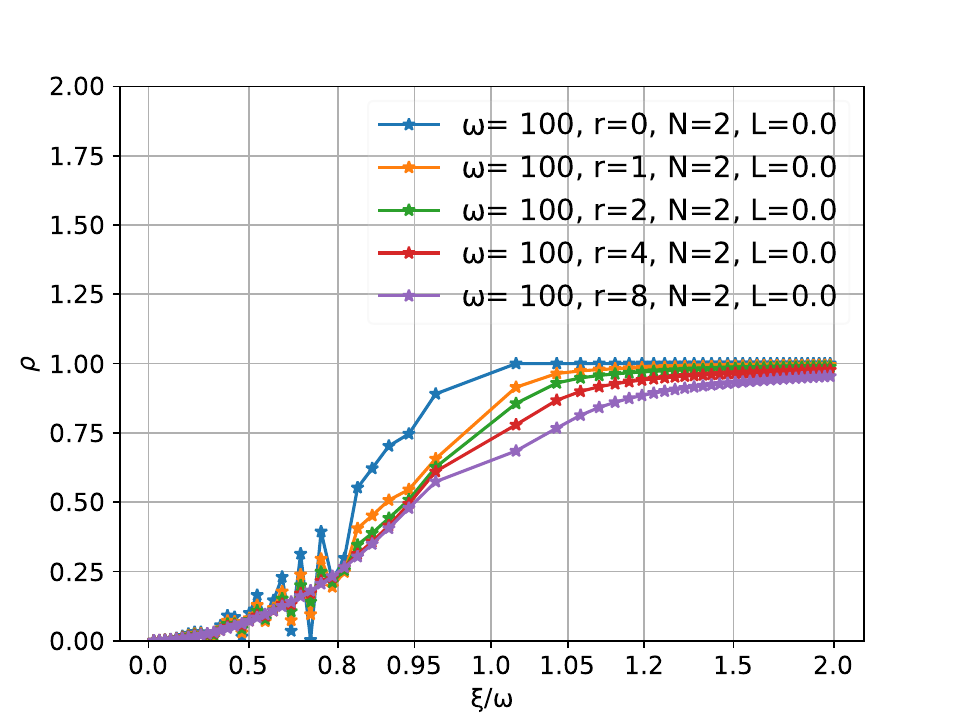}
  \includegraphics[width=0.33\textwidth,trim=10 5 45 37,clip]{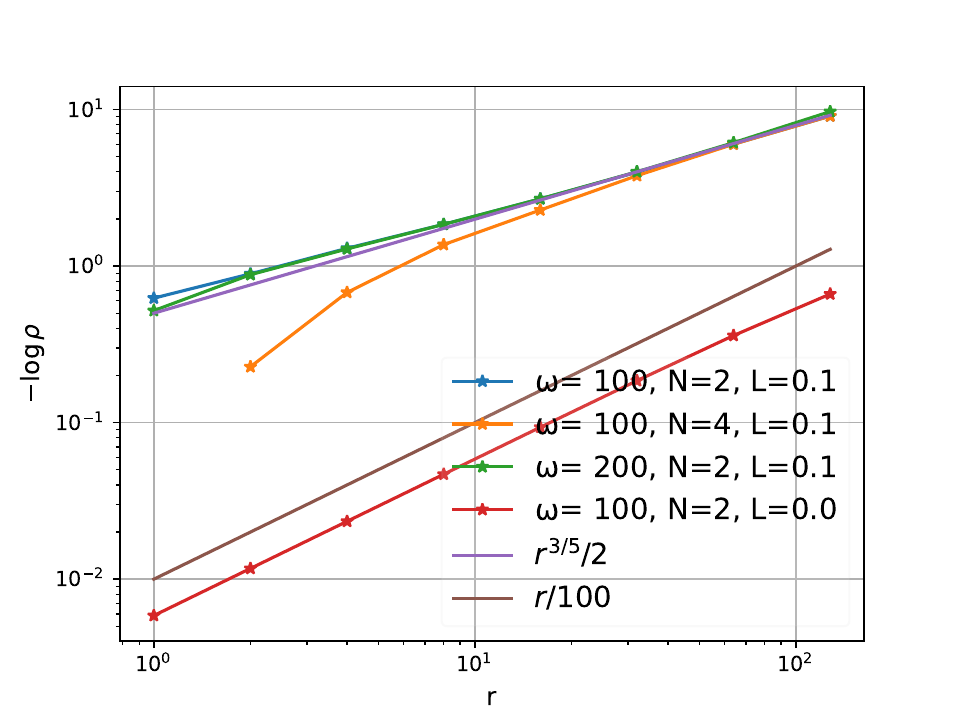}}
  \caption{Convergence factor dependence on $r$ for the waveguide with the operator
    $\Delta + \omega^2 - \I\omega r$}
  \label{fig1}
\end{figure}

We first measure how the convergence factor scales with the damping
coefficient $r$; see Fig.~\ref{fig1}.
We see that the damping helps not only convergence of the propagating
modes ($\xi/\omega<1$) but also the evanescent modes ($\xi/\omega>1$),
and that the benefit is more notable for the nonoverlapping method
($L=0$) and its evanescent modes. On the right, we see that the
convergence factor goes to zero exponentially when $r\to \infty$,
indicating very fast convergence.

\begin{figure}[t]
  \centering
  \mbox{\includegraphics[width=1.5in,trim=10 5 42 37,clip]{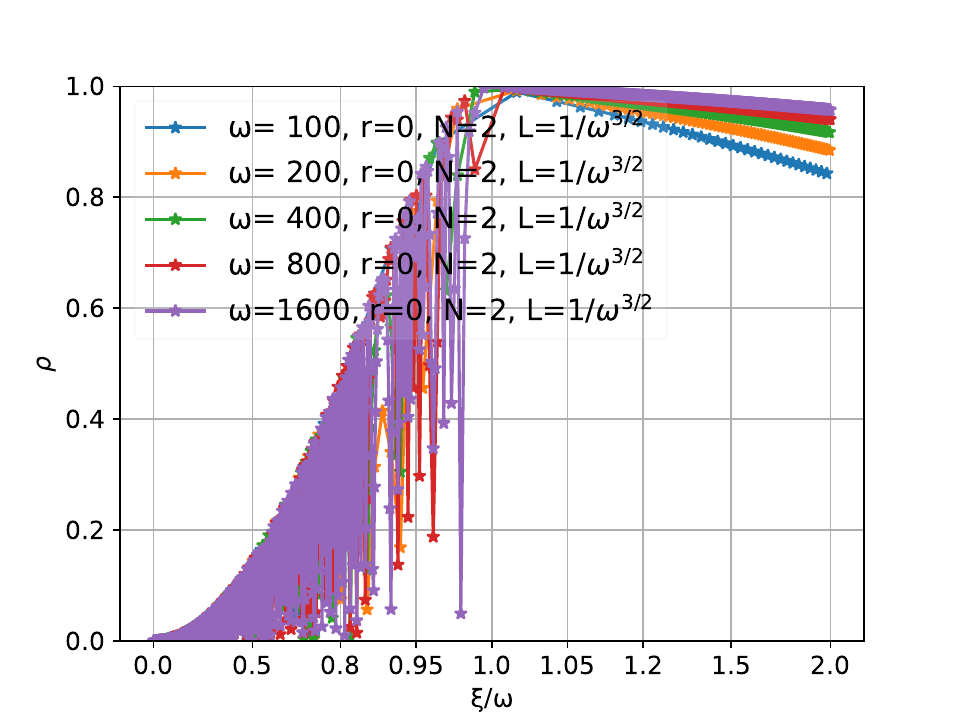}%
  \includegraphics[width=1.5in,trim=10 5 42 37,clip]{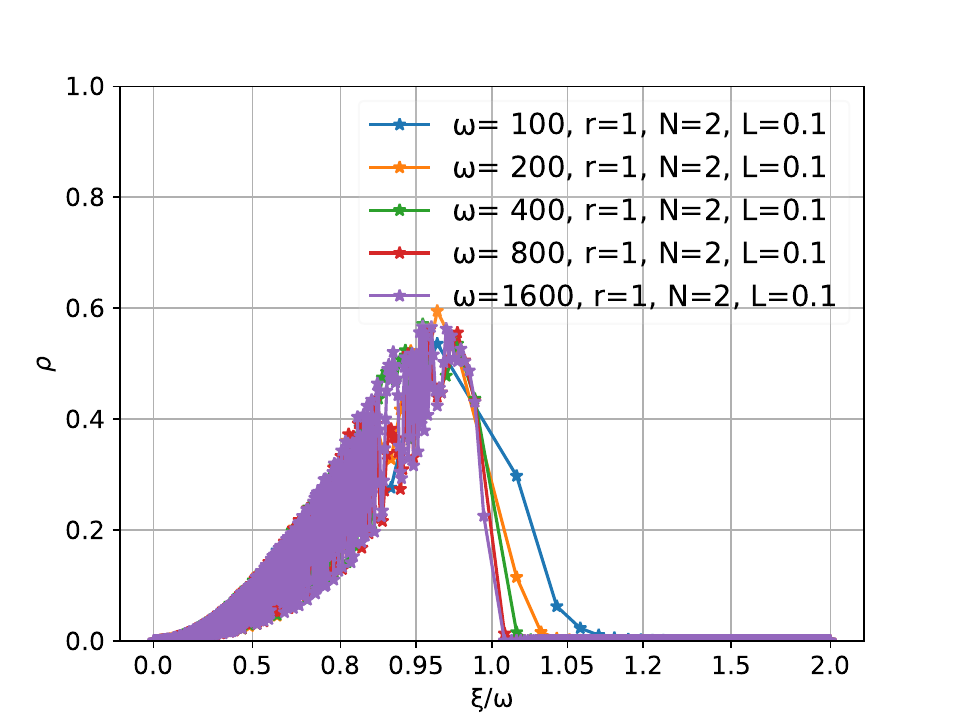}%
  \includegraphics[width=1.5in,trim=10 5 45 37,clip]{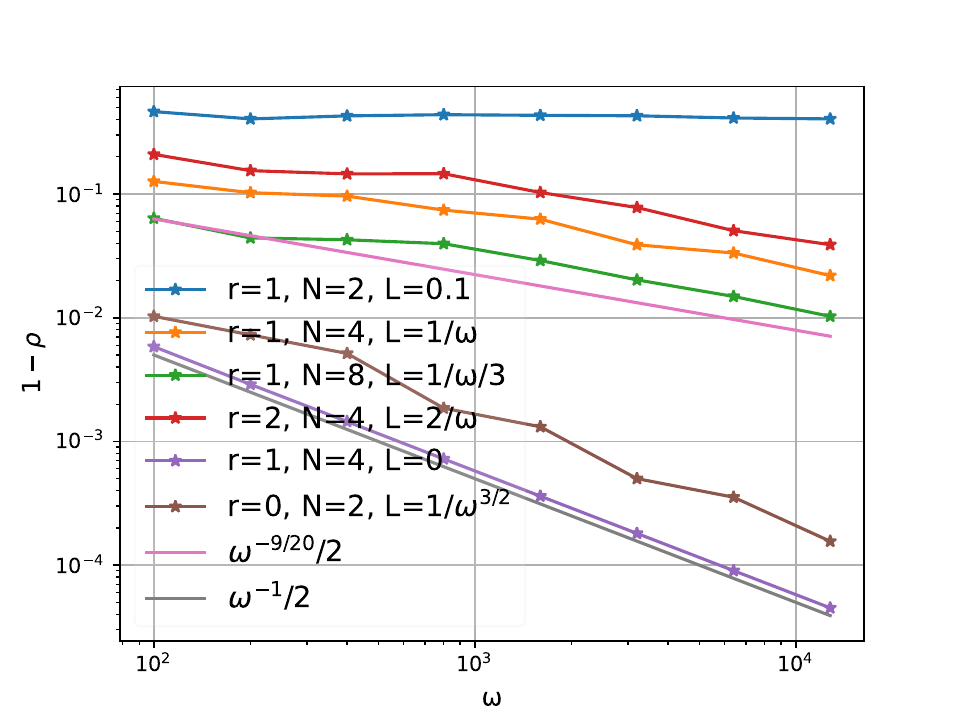}}
  \mbox{\includegraphics[width=1.5in,trim=10 5 42 37,clip]{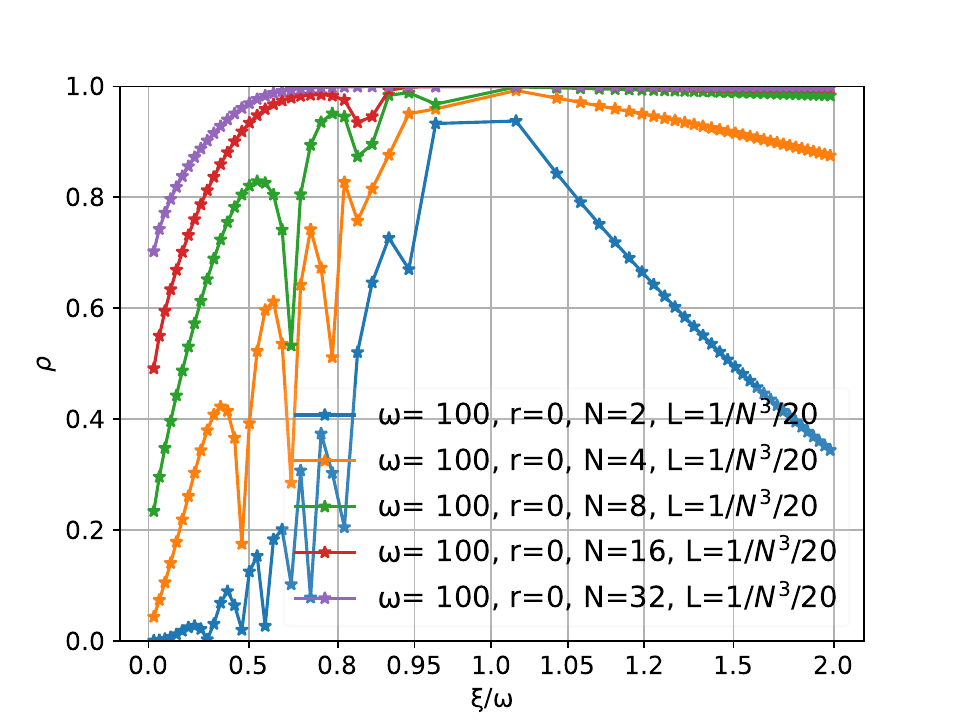}%
  \includegraphics[width=1.5in,trim=10 5 42 37,clip]{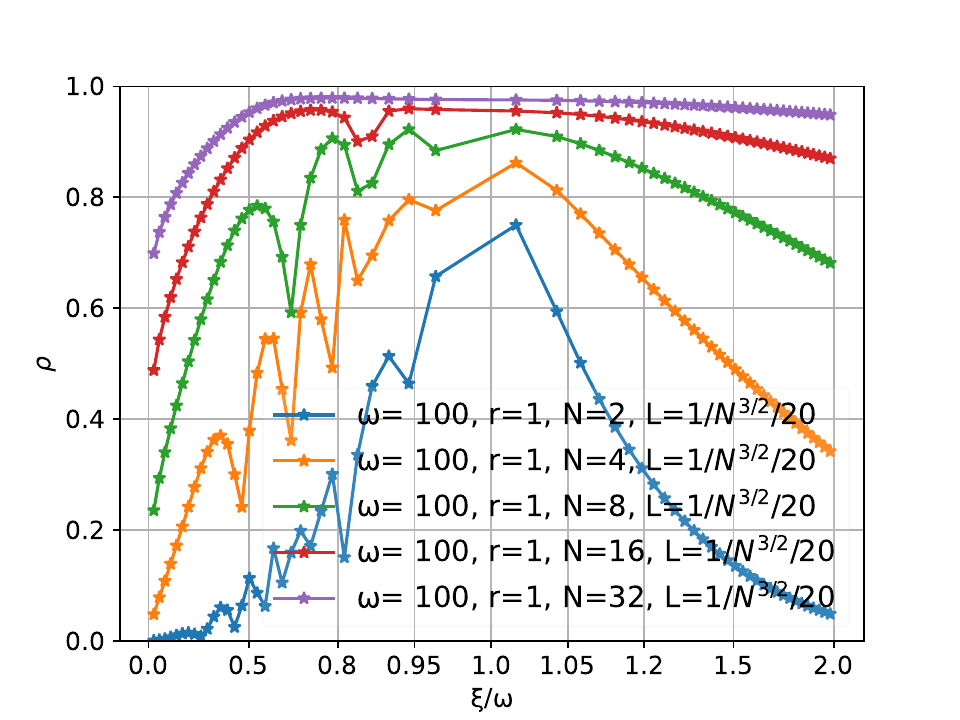}%
  \includegraphics[width=1.5in,trim=10 5 45 37,clip]{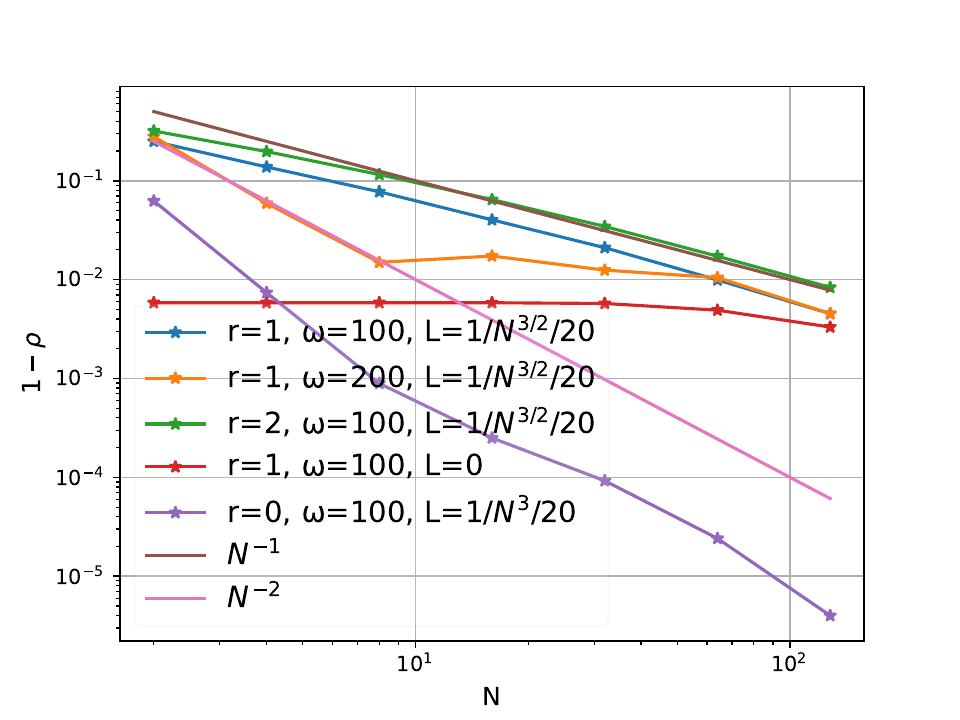}}
  \mbox{\includegraphics[width=1.5in,trim=10 5 42 37,clip]{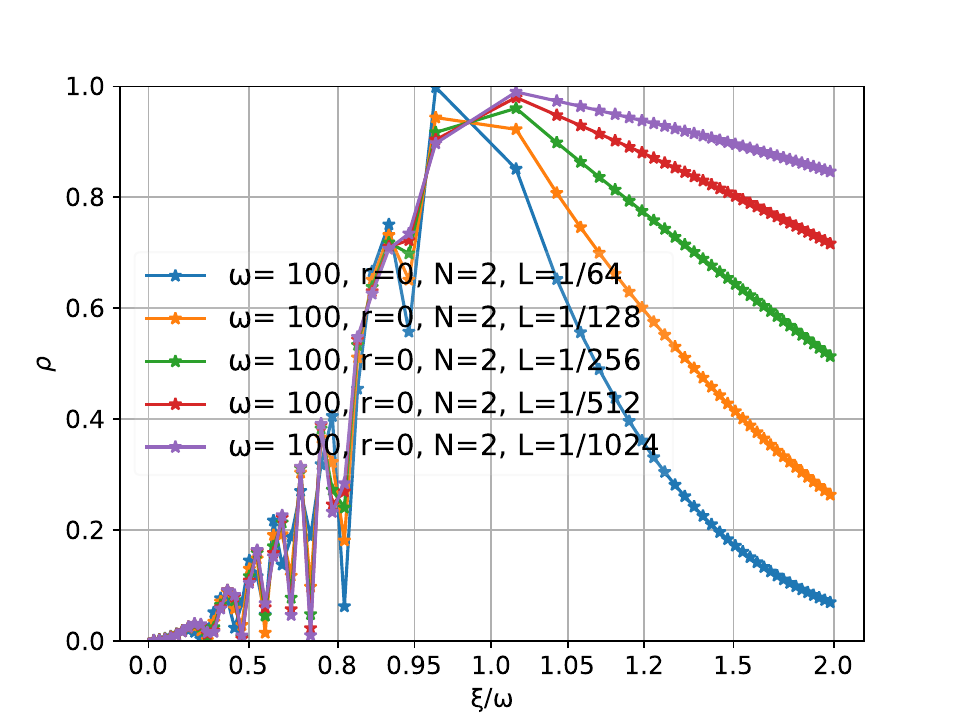}%
  \includegraphics[width=1.5in,trim=10 5 42 37,clip]{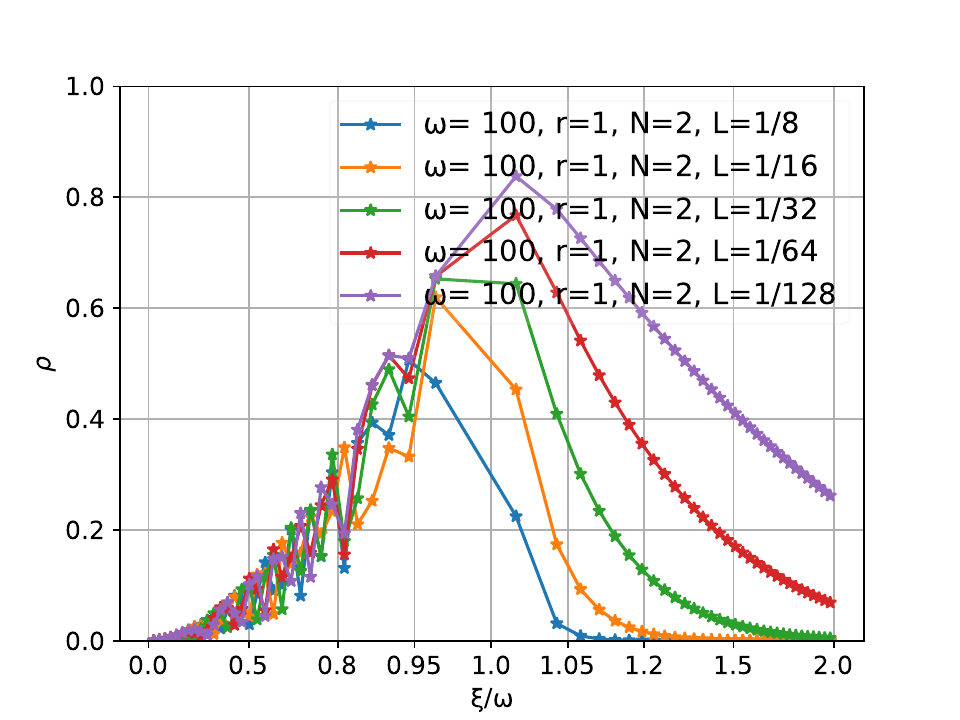}%
  \includegraphics[width=1.5in,trim=10 5 45 37,clip]{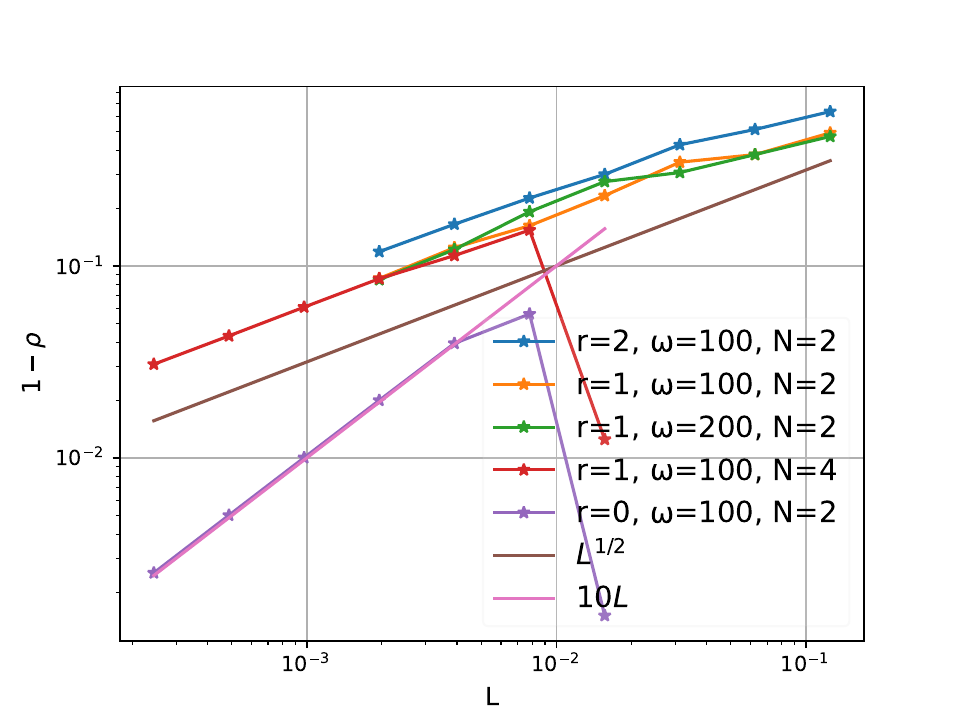}}%
  \caption{Convergence factor dependence on $\omega$ (top), number of
    subdomains $N$ (middle) and overlap $L$ (bottom) for the waveguide
    with the operator $\Delta + \omega^2 - \I\omega r$}
  \label{fig24}
\end{figure}

Next, we study the scalings with respect to the wavenumber $\omega$,
number of subdomains $N$ and the overlap width $L$. On many
subdomains, the overlap width $L$ has to be sufficiently small with
respect to the wavenumber $\omega$ and the number of subdomains $N$ to
ensure convergence without physical damping
\cite{Gander_Zhang_2022}. Actually, for each given $\omega$ and $N$,
there would be an optimal overlap width, but we will not seek the
optimal $L$ in this note and just use an ad hoc choice.
Fig.~\ref{fig24}
show that the waveguide problem with first order damping has similar
scalings to the free space problem without damping in
\cite{Gander_Zhang_2022}, while the waveguide problem without damping
has much different scalings and is more difficult to solve by the
Schwarz method. In particular, the top row in Fig.~\ref{fig24} shows
that the Schwarz method on two fixed subdomains is robust in the
wavenumber with damping, but not on many subdomains or without
damping.

\subsection{Helmholtz operator $(1+ \I\omega \gamma)\Delta + \omega^2$ in the waveguide}

\begin{figure}[t]
  \centering
  \includegraphics[width=1.5in,trim=10 5 20 35,clip]{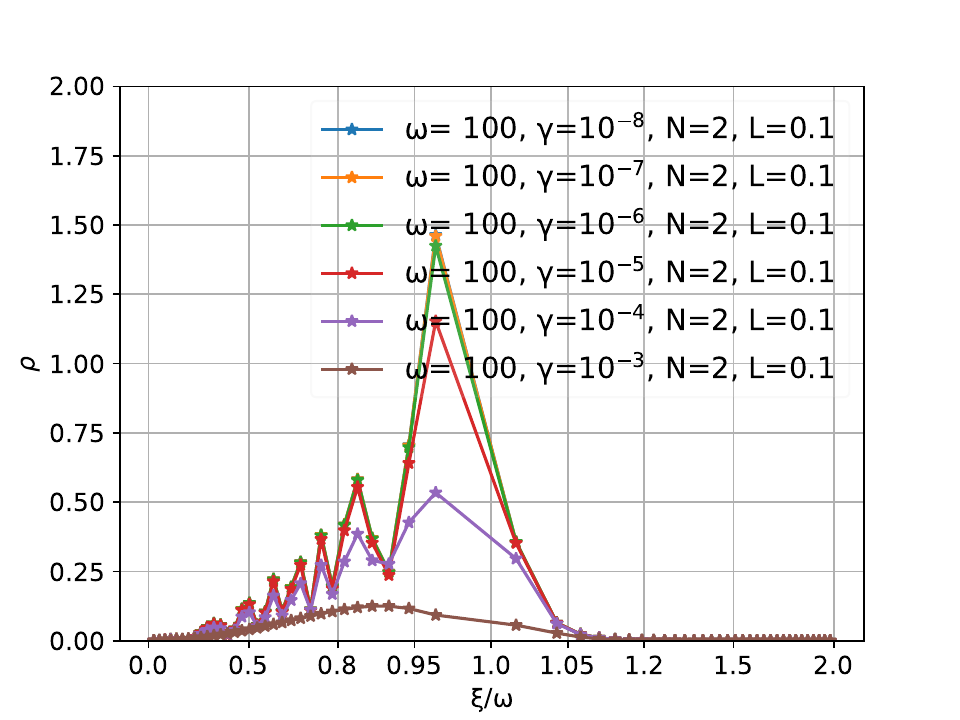}%
  \includegraphics[width=1.5in,trim=10 5 20 35,clip]{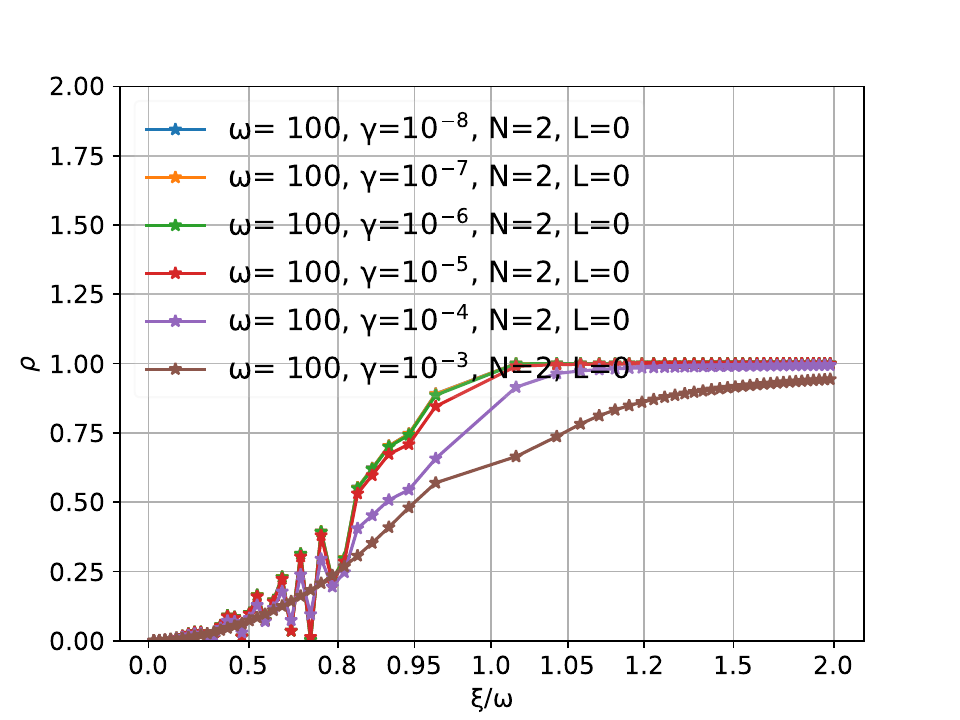}
  \includegraphics[width=1.5in,trim=10 5 20 35,clip]{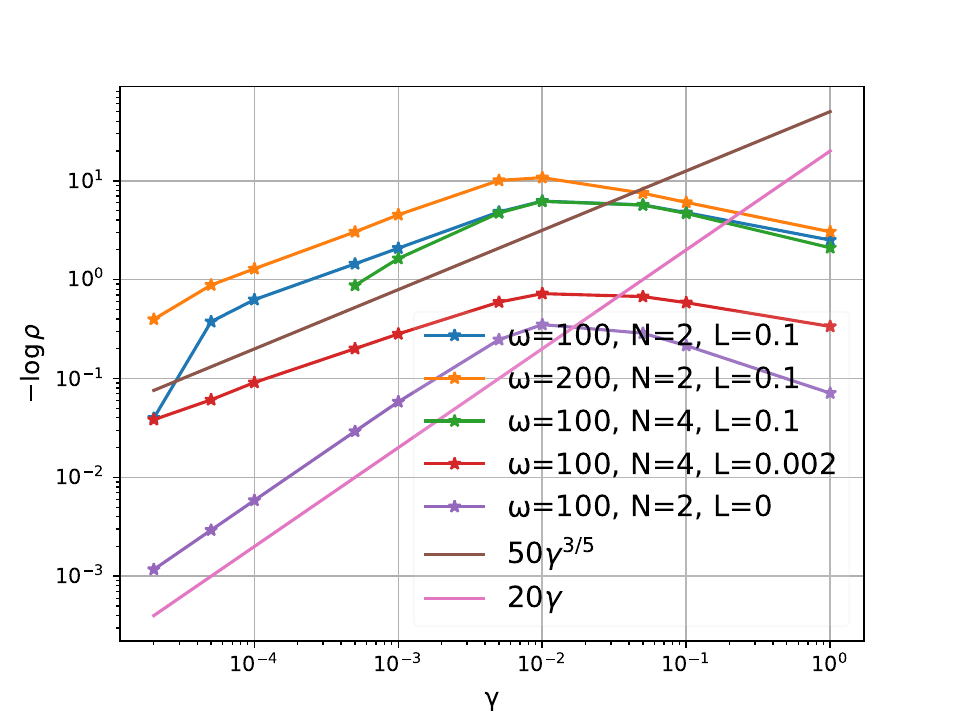}
  \caption{Convergence factor dependence on $\gamma$ for the waveguide
    with the operator $(1+ \I\omega \gamma)\Delta + \omega^2$}
  \label{fig5}
\end{figure}

As we mentioned, the Helmholtz equation $(1+ \I\omega \gamma)\Delta u
+ \omega^2u=f$ may be normalized to $\Delta u + u\omega^2/(1+ \I\omega
\gamma) = f/(1+ \I\omega \gamma)$. The 0th order coefficient is then
\[
  \omega^2/(1+ \I\omega \gamma) = \frac{\omega^2}{1+\omega^2\gamma^2}(1-\I\omega \gamma)
  \approx \begin{cases} \omega^2-\I\omega^3\gamma & \text{if }\omega \gamma\ll 1\\
    c\omega^2(1-\I) & \text{if }\omega\gamma\approx 1\text{ and }\frac{1}{1+\omega^2\gamma^2}=c\\
  \gamma^{-2}-\I\omega\gamma^{-1} & \text{if }\omega\gamma\gg 1\end{cases},
\]
and, with no approximation, the ratio of the imaginary part to the
real part is $-\omega\gamma$.  Fig.~\ref{fig5}
shows the scaling of the convergence factor with $\gamma$. Note that
$\gamma=0$ here is the same as the no damping case $r=0$ in
Sect.~\ref{sec2.1}, so $\gamma=0$ is not shown here. But compared to
Fig.~\ref{fig1} with $\omega=100$ the small dampings $\gamma=10^{-8}$,
$10^{-7}$, $10^{-6}$ almost coincide with the no damping case, while
$\gamma=10^{-4}$ roughly mimics $r=\omega^2\gamma=1$.

\begin{figure}[t]
  \centering
  \mbox{\includegraphics[width=1.5in,trim=10 5 42 37,clip]{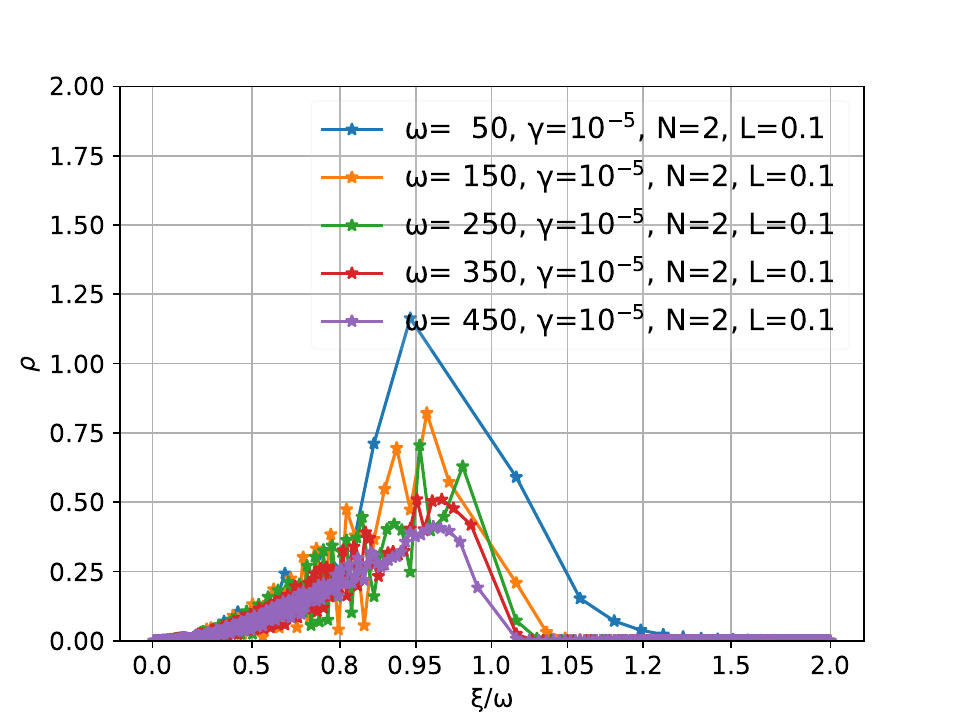}%
  \includegraphics[width=1.5in,trim=10 5 42 37,clip]{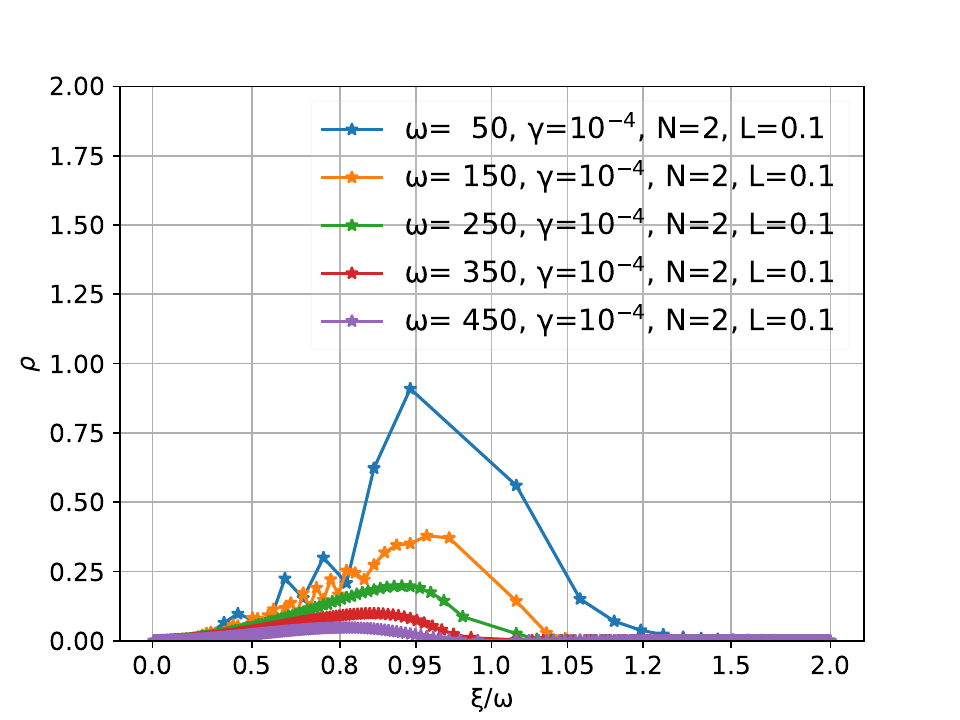}
  \includegraphics[width=1.5in,trim=10 5 45 37,clip]{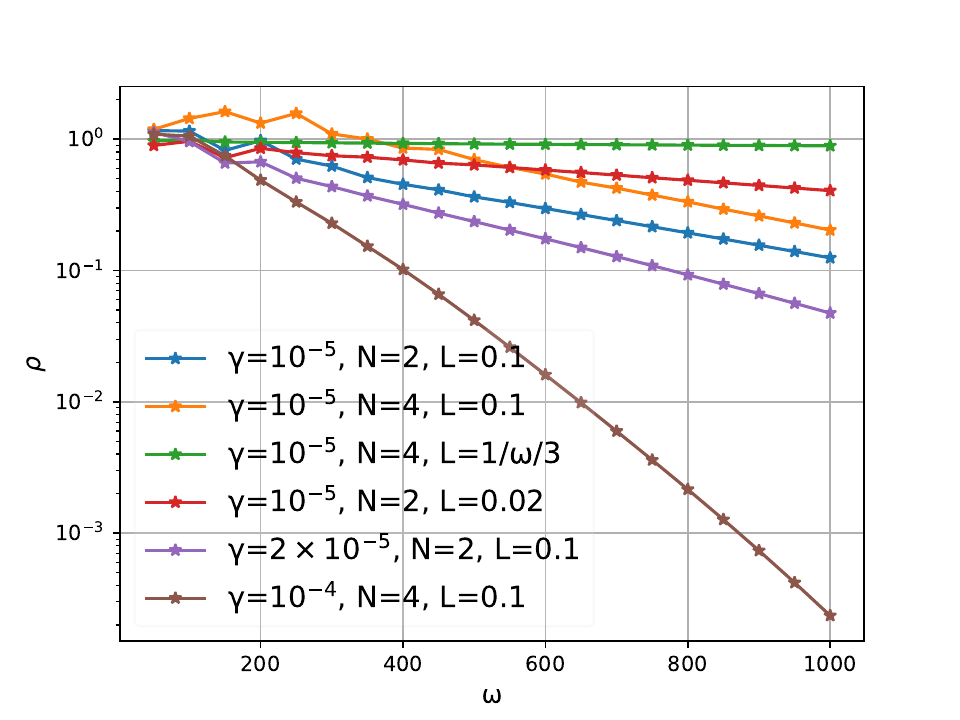}}
  \mbox{\includegraphics[width=1.5in,trim=10 5 42 37,clip]{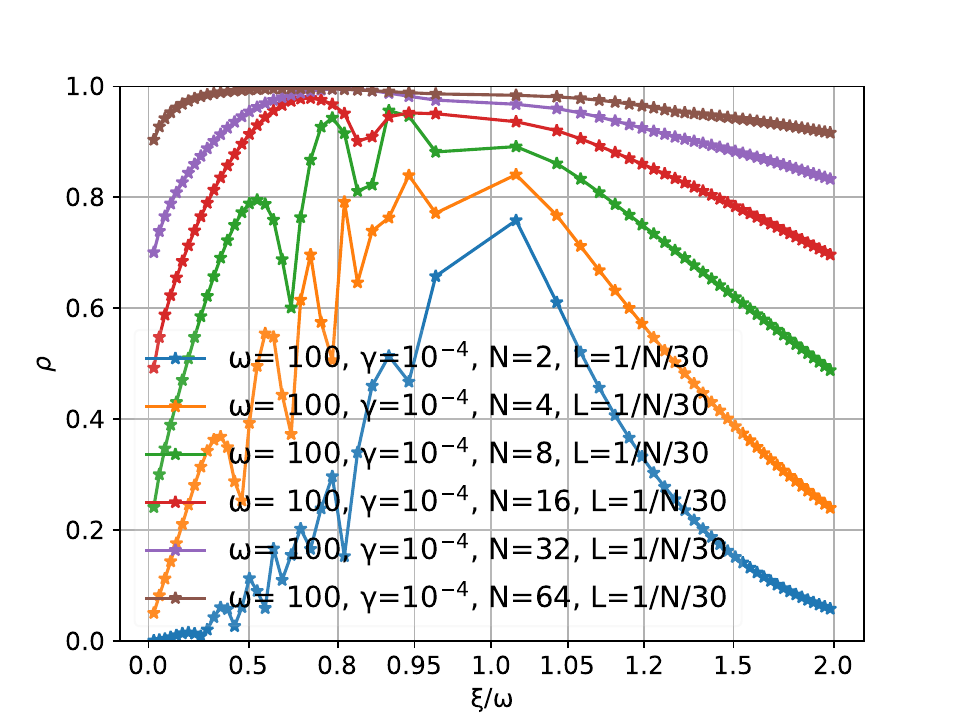}%
  \includegraphics[width=1.5in,trim=10 5 42 37,clip]{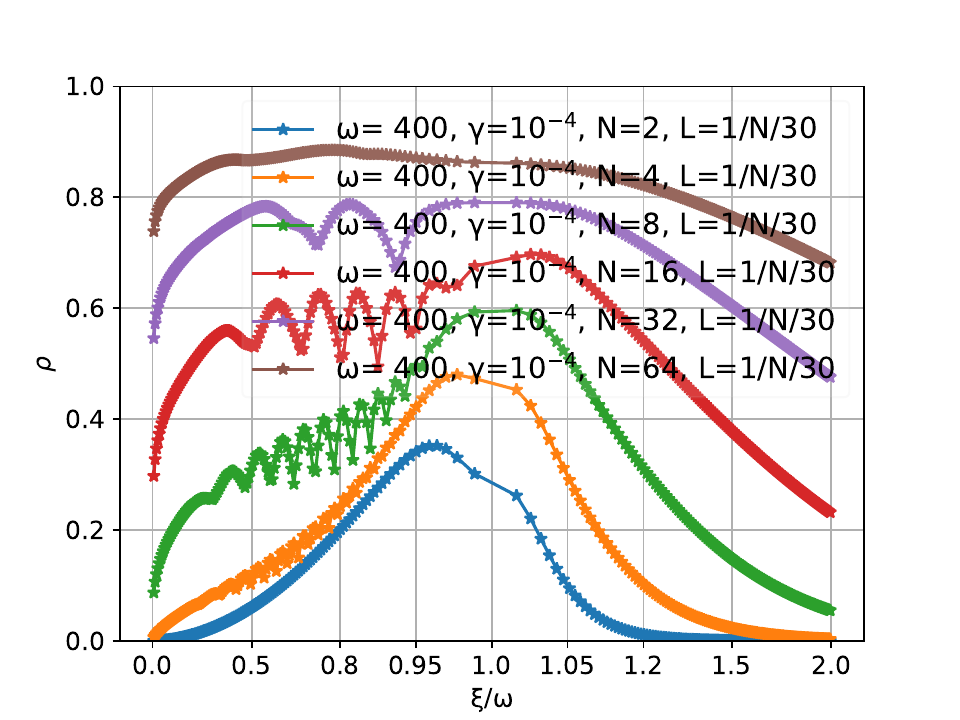}
  \includegraphics[width=1.5in,trim=10 5 45 37,clip]{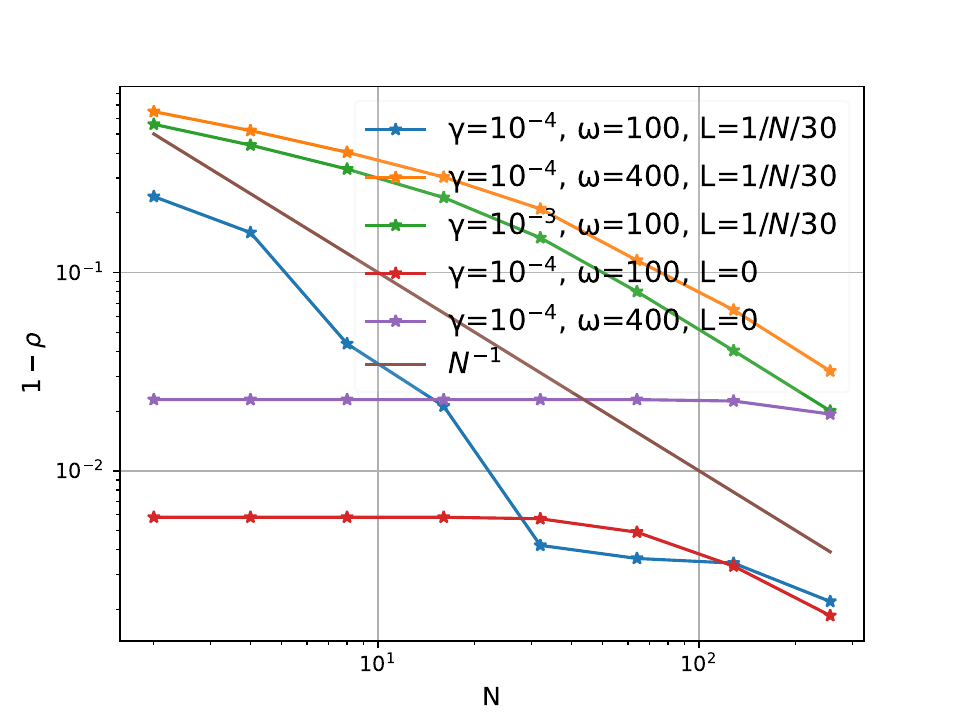}}
  \mbox{\includegraphics[width=1.5in,trim=10 5 42 37,clip]{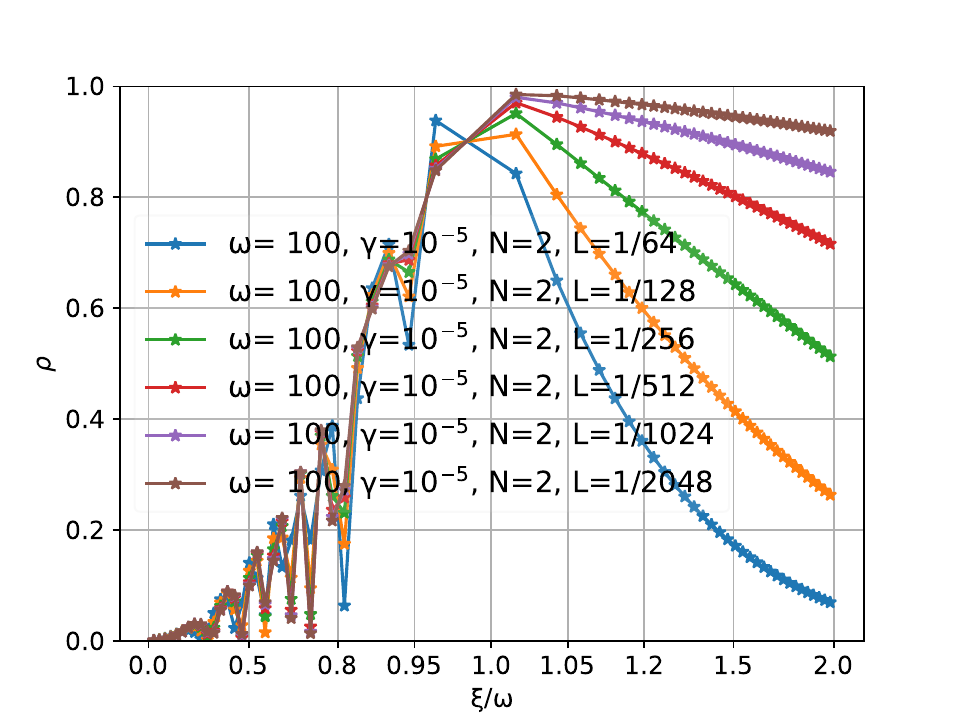}%
  \includegraphics[width=1.5in,trim=10 5 42 37,clip]{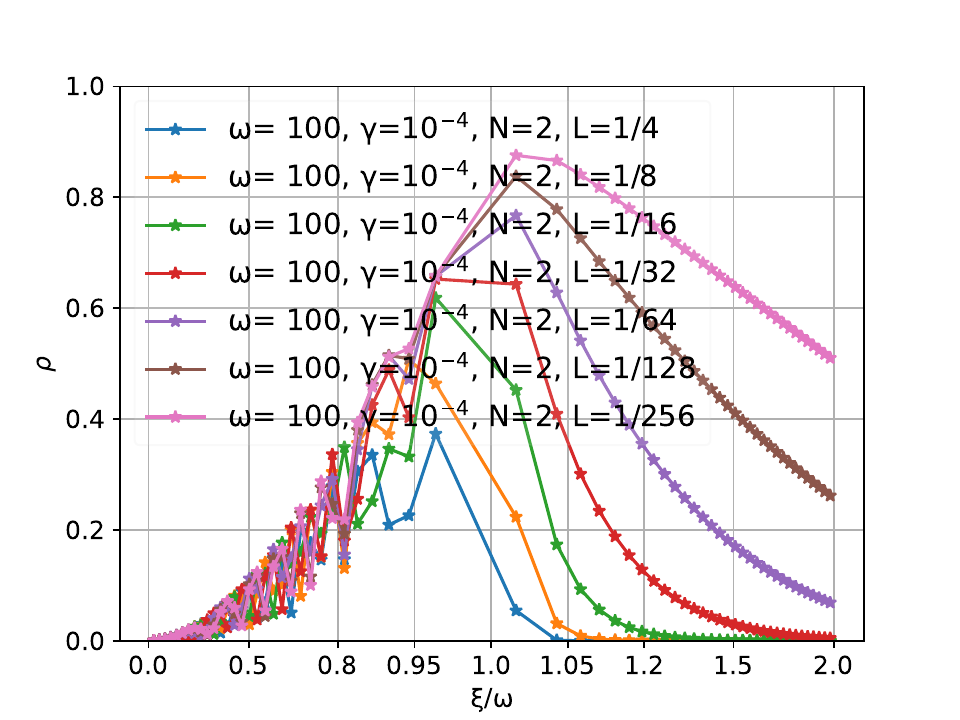}
  \includegraphics[width=1.5in,trim=10 5 45 37,clip]{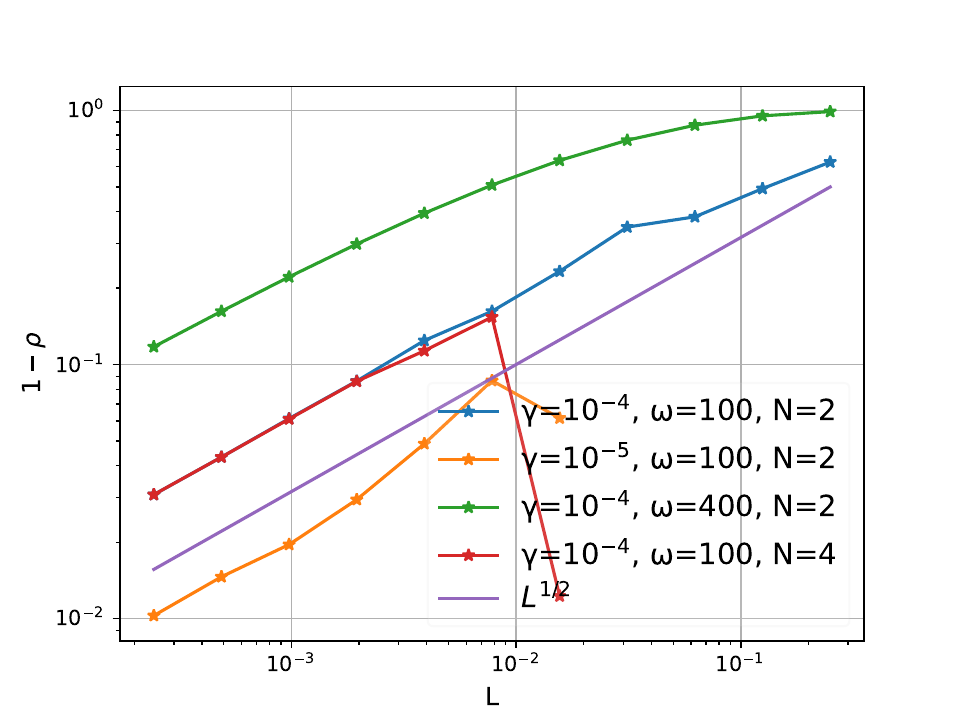}}
  \caption{Convergence factor dependence on $\omega$ (top), number of
    subdomains $N$ (middle) and overlap $L$ (bottom) for the waveguide
    with the operator $(1+ \I\omega \gamma)\Delta + \omega^2$}
  \label{fig68}
\end{figure}

However, this does not imply that the viscoelastic damping in $(1+
\I\omega \gamma)\Delta + \omega^2$ is the same as the first order
damping in $\Delta + \omega^2 - \I\omega r$: the latter has the
imaginary/real ratio $-r/\omega$. So we can expect the two types of
damping to differ significantly with large $\omega$. This is clearly
visible in the first line of Fig.~\ref{fig68}
which is drastically different from the first line in
Fig.~\ref{fig24}. The scaling with number of subdomains is shown in
Fig.~\ref{fig68} (middle). For moderate wavenumber, sufficiently small
overlap is still required for convergence on many subdomains. There is
an optimal overlap width $L$ for given $\omega$, $\gamma$ and $N$,
which is not investigated here. But with the ad hoc choice of
$L=O(N^{-1})$ we get the typical scaling $\rho=1-O(N^{-1})$ for
one-level parallel Schwarz methods applied to Laplace problems.

In the bottom row of Fig.~\ref{fig68}, we show the dependence on the
overlap width $L$. The scaling with small overlap is the same as
Fig.~\ref{fig24} (bottom), also the same requirement of smaller overlap
for more subdomains. But here larger wavenumbers can lead to faster
convergence and allow using larger overlap.

\section{Cavity Problem}\label{sec3}

In the closed cavity problem, $u=0$ on all the boundary. For
wellposedness of the problem without damping, we assume the squared
wavenumber $\omega^2$ is not an eigenvalue of the negative Laplace
operator.

\subsection{Helmholtz operator $\Delta + \omega^2 - \I\omega r$ in the cavity}
\label{sec3.1}

\begin{figure}[t]
  \centering
  \includegraphics[width=1.5in,trim=10 5 42 37,clip]{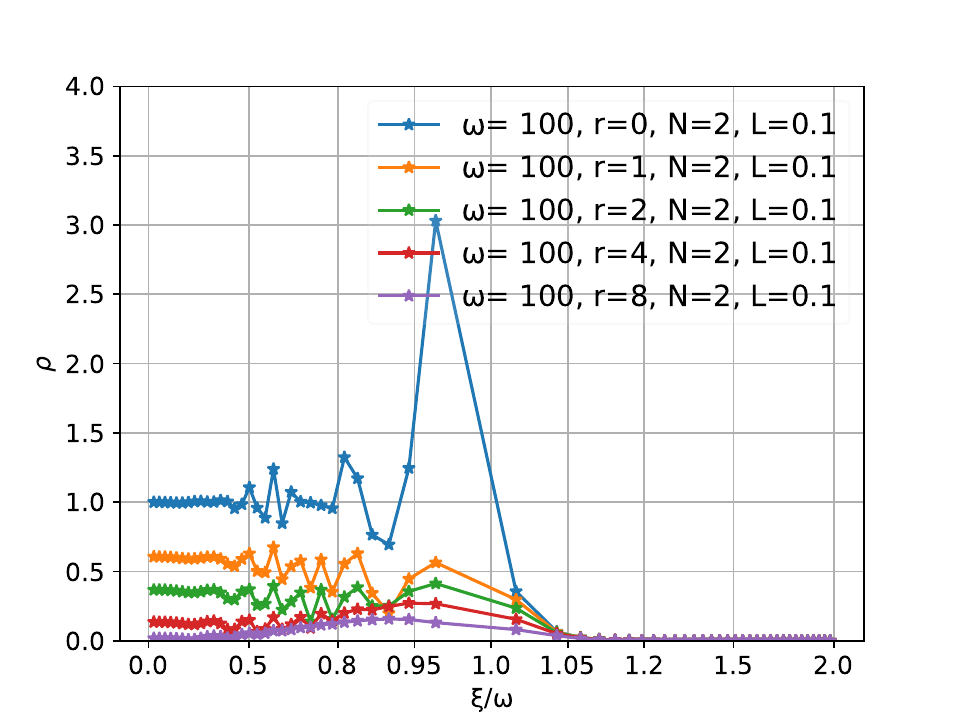}%
  \includegraphics[width=1.5in,trim=10 5 42 37,clip]{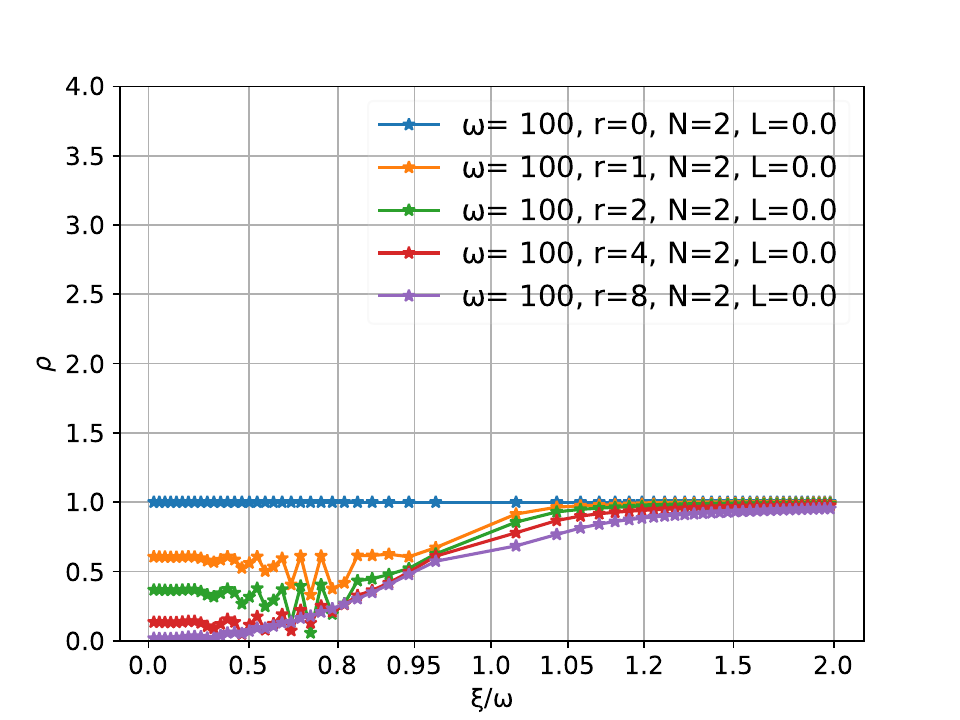}
  \includegraphics[width=1.5in,trim=10 5 45 37,clip]{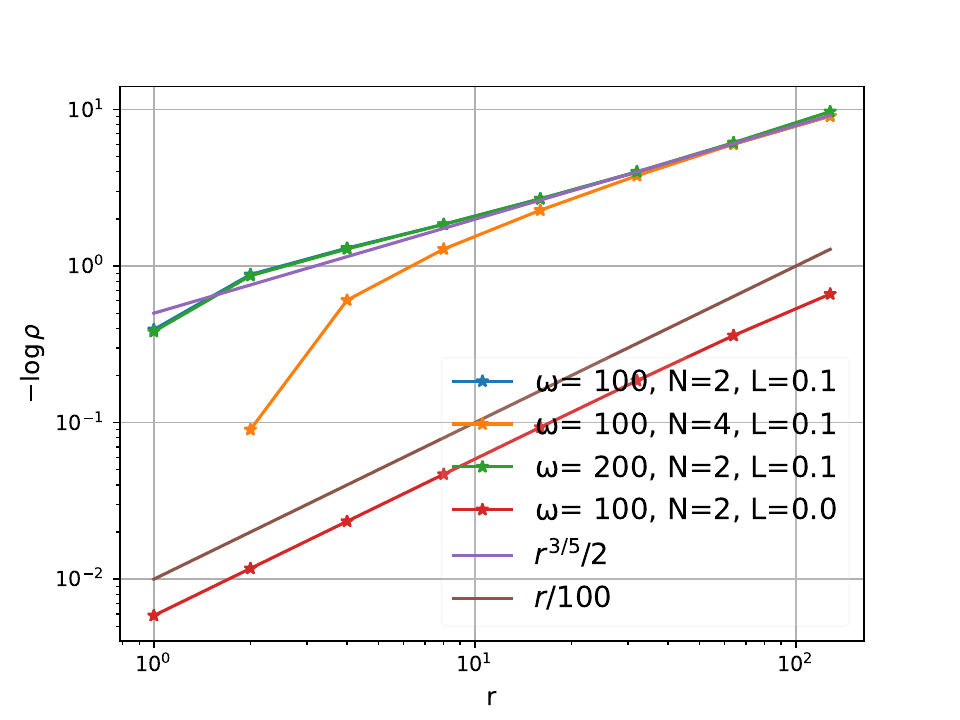}
  \caption{Convergence factor dependence on $r$ for the cavity with
    the operator $\Delta + \omega^2 - \I\omega r$}
  \label{fig9}
\end{figure}

The influence of the first order damping coefficient $r$ is shown in
Fig.~\ref{fig9}.
Compared to Fig.~\ref{fig1} for the waveguide problem, here for the
cavity problem the low space frequency has a larger convergence
factor. But this difference quickly diminishes as $r$ increases, and
the scaling for $r$ large remains the same, so even the very hard
closed cavity problem becomes easy with first order damping.

\begin{figure}[t]
  \centering
  \mbox{\includegraphics[width=1.5in,trim=10 5 42 37,clip]{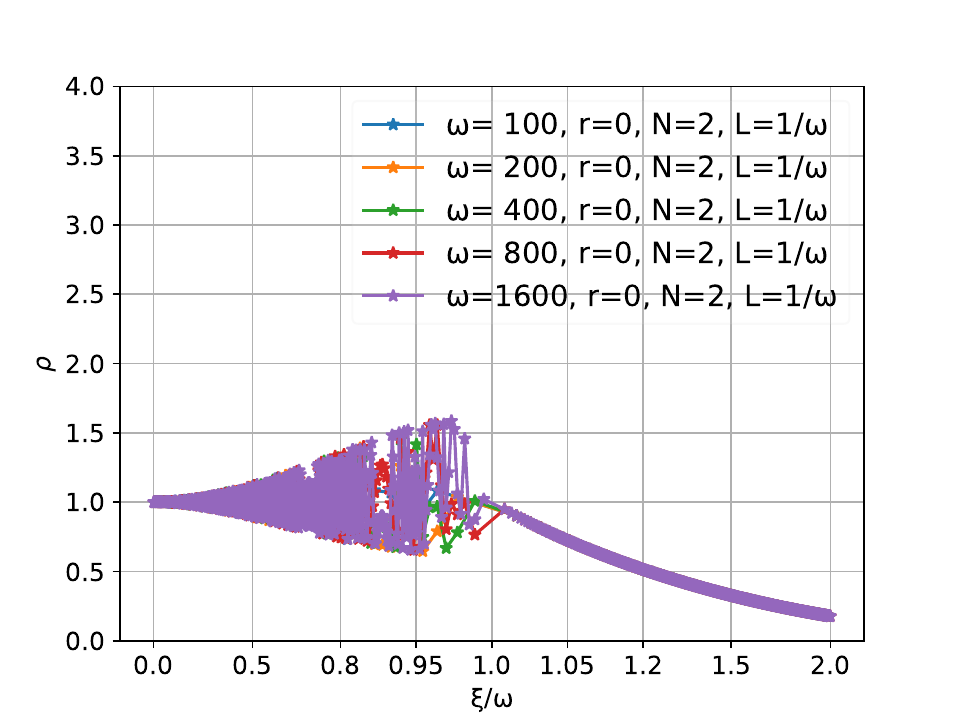}%
  \includegraphics[width=1.5in,trim=10 5 42 37,clip]{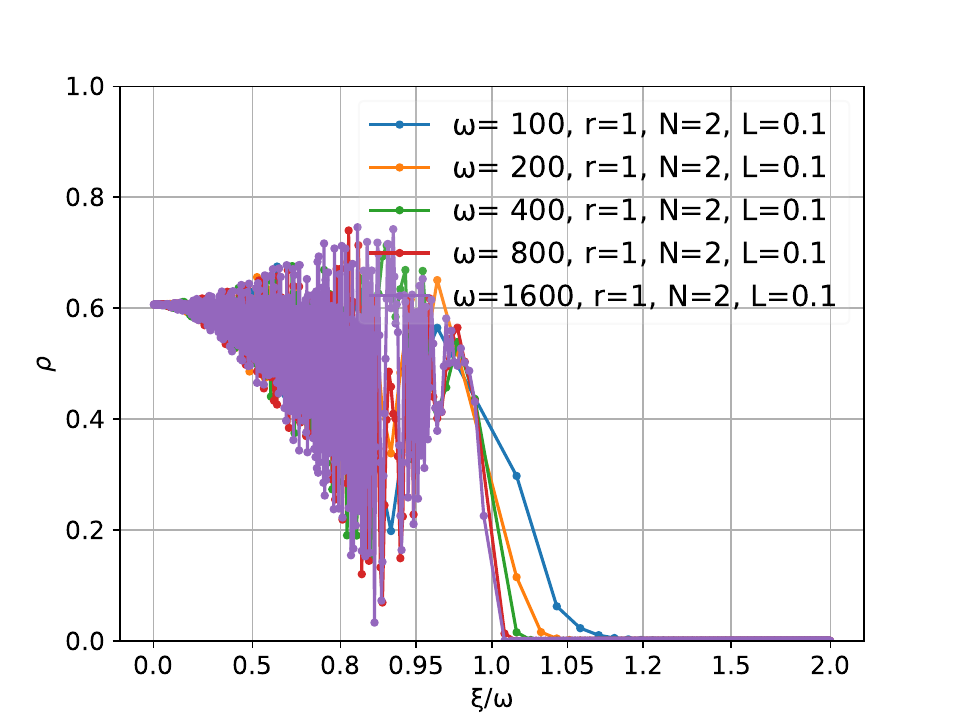}%
  \includegraphics[width=1.5in,trim=10 5 45 37,clip]{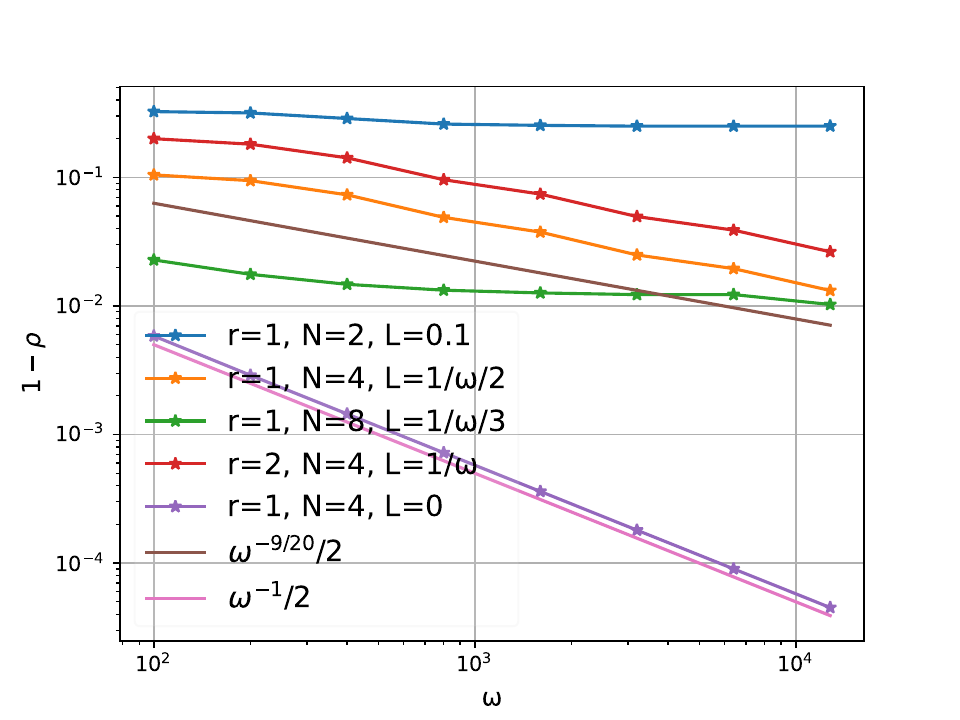}}
  \mbox{\includegraphics[width=1.5in,trim=10 5 42 37,clip]{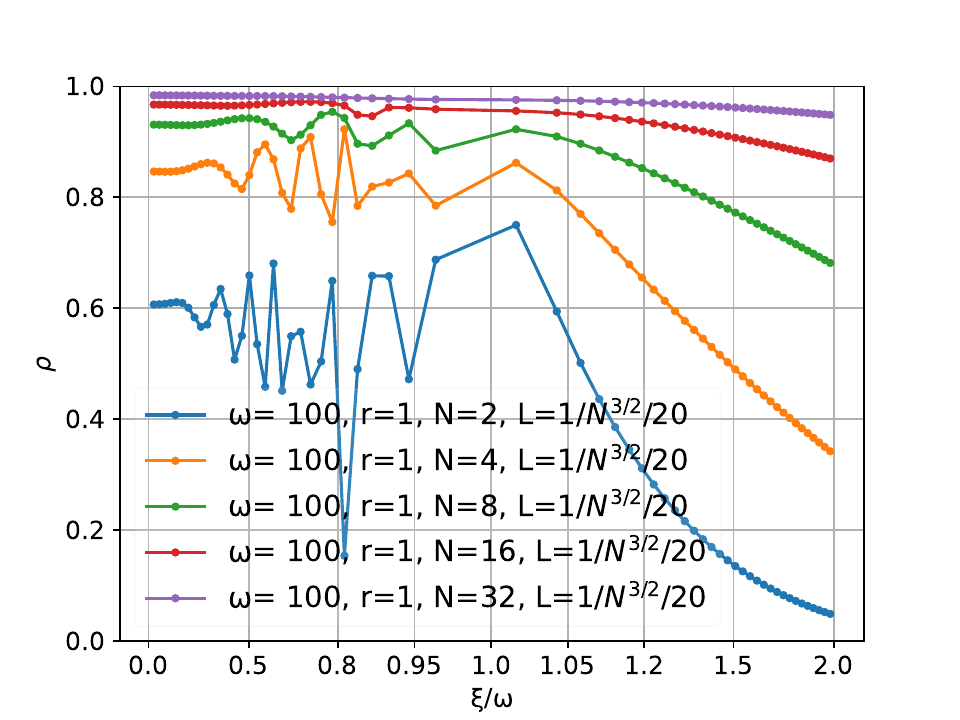}%
  \includegraphics[width=1.5in,trim=10 5 42 37,clip]{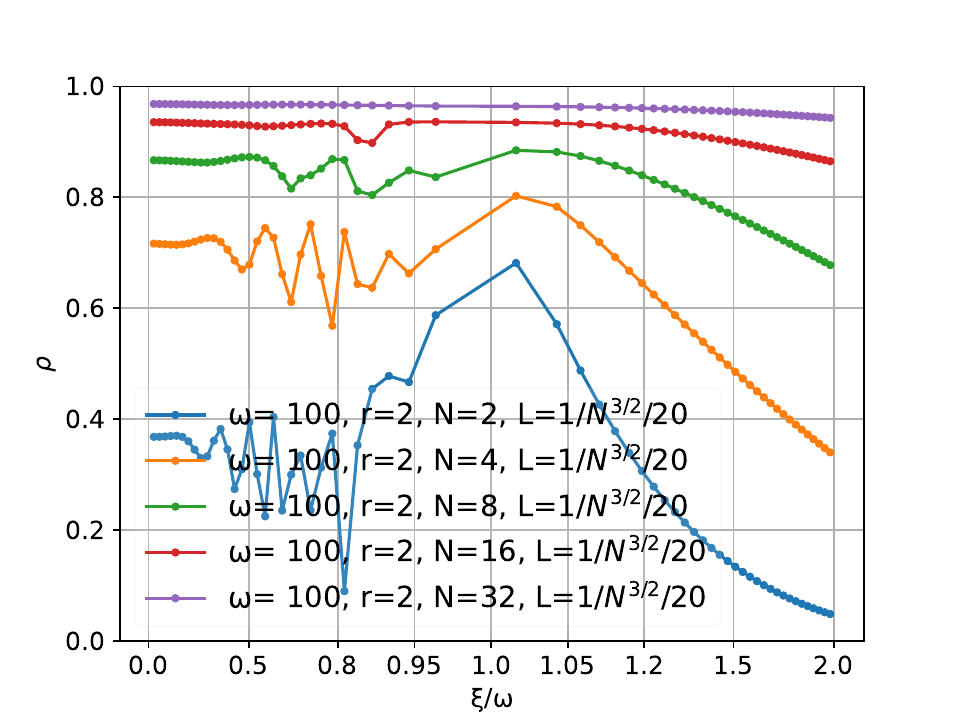}%
  \includegraphics[width=1.5in,trim=10 5 45 37,clip]{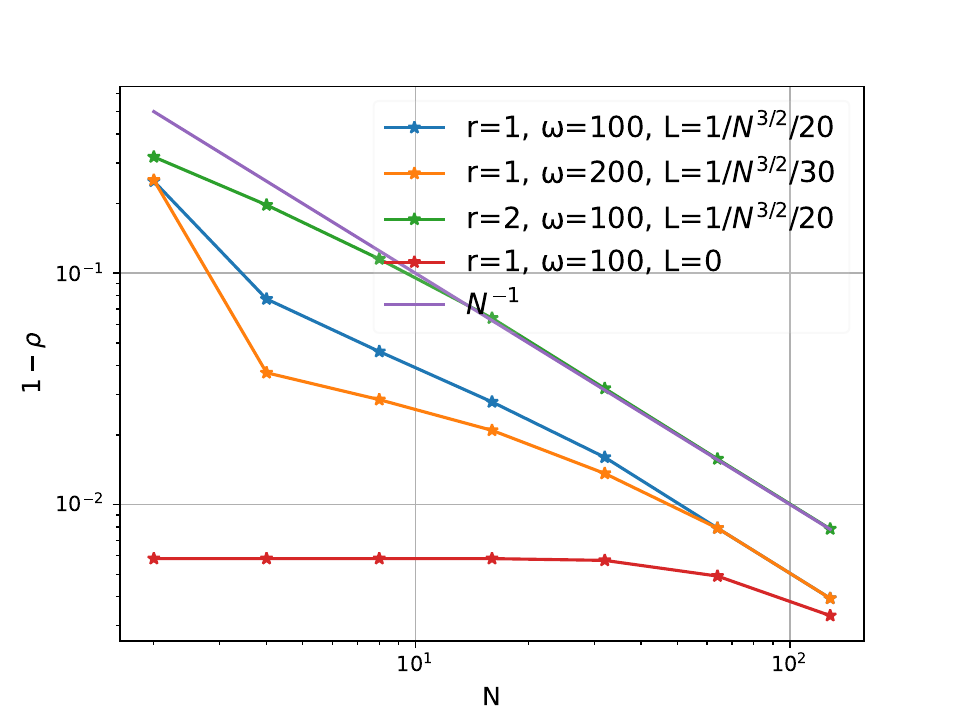}}
  \mbox{\includegraphics[width=1.5in,trim=10 5 42 37,clip]{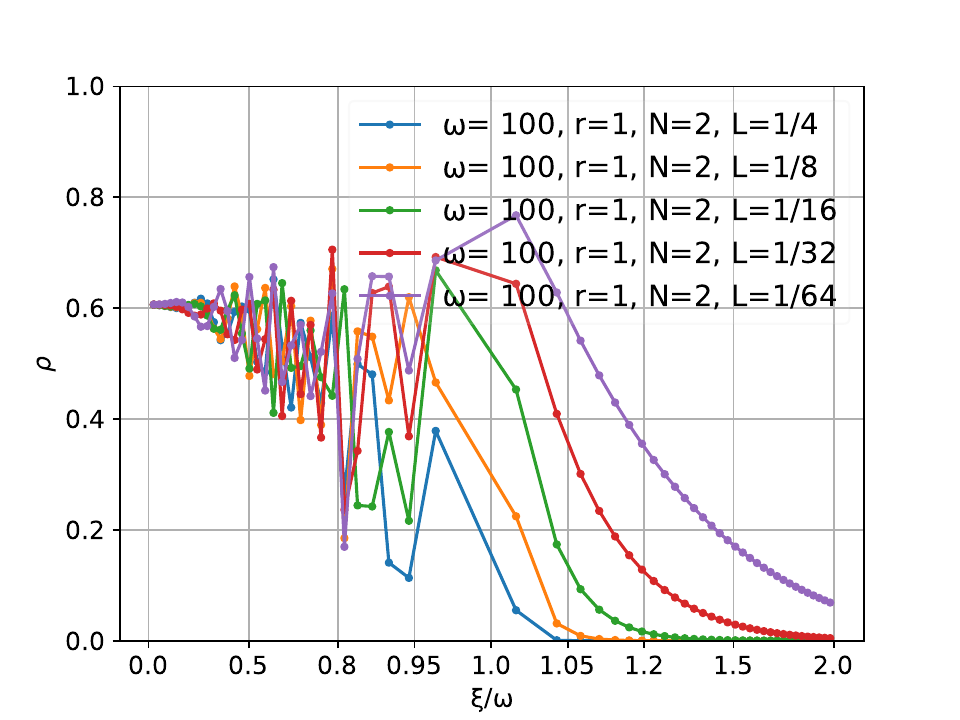}%
  \includegraphics[width=1.5in,trim=10 5 42 37,clip]{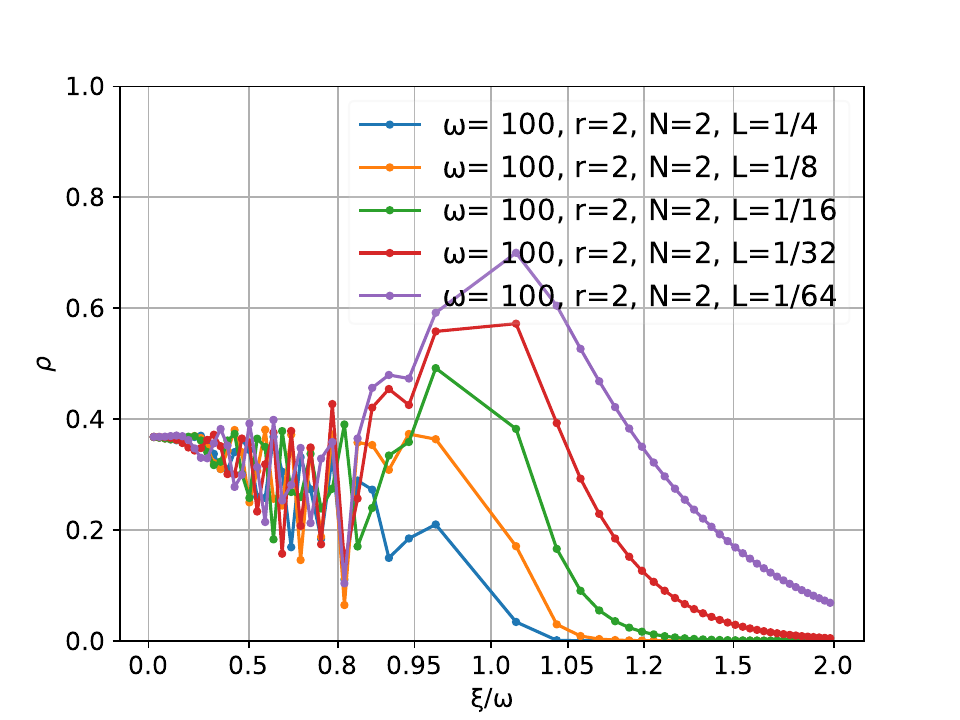}%
  \includegraphics[width=1.5in,trim=10 5 45 37,clip]{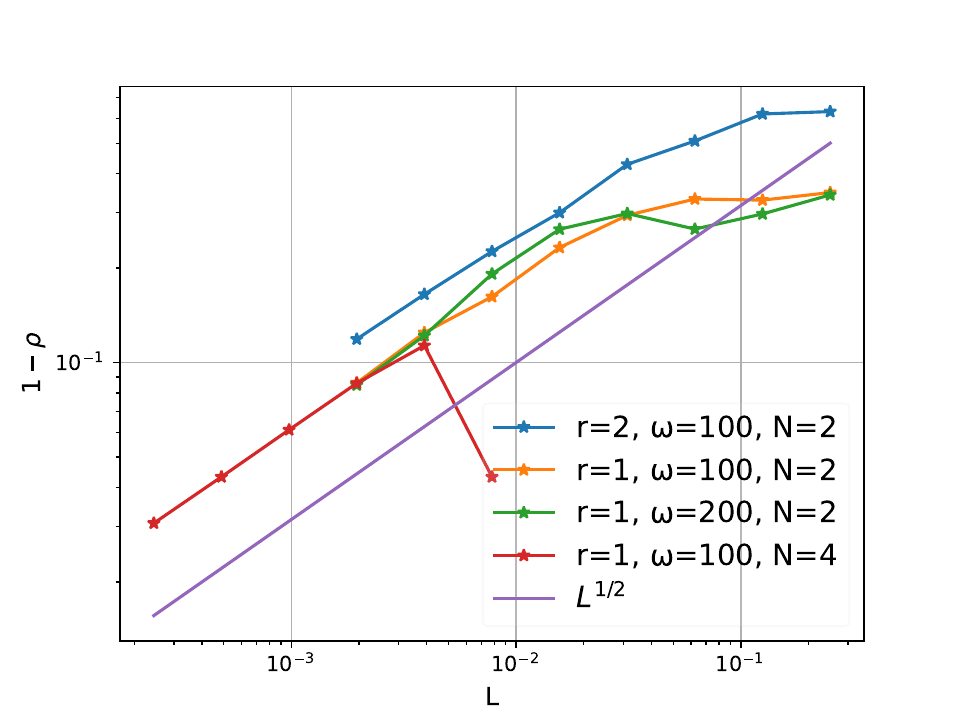}}%
  \caption{Convergence factor dependence on $\omega$ (top), number of
    subdomains $N$ (middle) and overlap $L$ (bottom) for the cavity
    with the operator $\Delta + \omega^2 - \I\omega r$}
  \label{fig1012}
\end{figure}

We show in Fig.~\ref{fig1012}
the corresponding dependence on the other parameters.  We see that
without first order damping the Schwarz method for the closed cavity
problem can not converge at all; see Fig.~\ref{fig1012} (top
left).  With first order damping, convergence is robust in $\omega$ for
two subdomains as in the wave guide case, and deteriorates similarly for more
subdomains. {From the green line ($r=1$, $N=8$, $L=1/\omega/3$), we can see that, as
  $\omega$ gradually increases, the convergence here is first slower than, and then close to, the
  convergence in the wave guide case.}

\subsection{Helmholtz operator $(1+\I\omega \gamma )\Delta  + \omega^2$ in the cavity}

\begin{figure}[t]
  \centering
  \mbox{\includegraphics[width=1.5in,trim=10 5 42 37,clip]{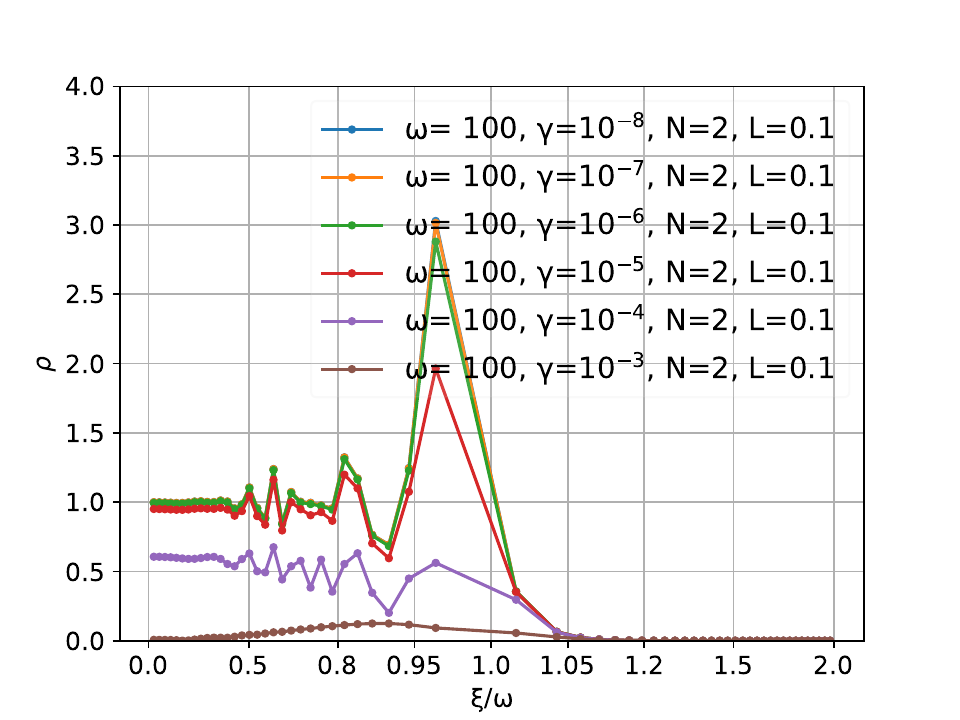}%
  \includegraphics[width=1.5in,trim=10 5 42 37,clip]{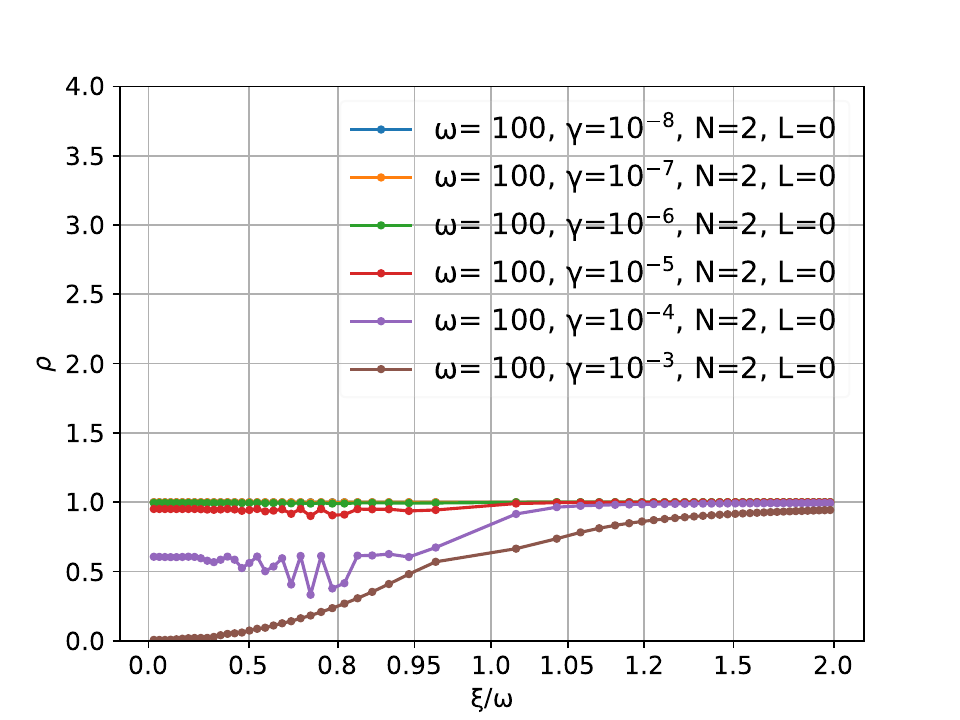}
  \includegraphics[width=1.5in,trim=10 5 45 37,clip]{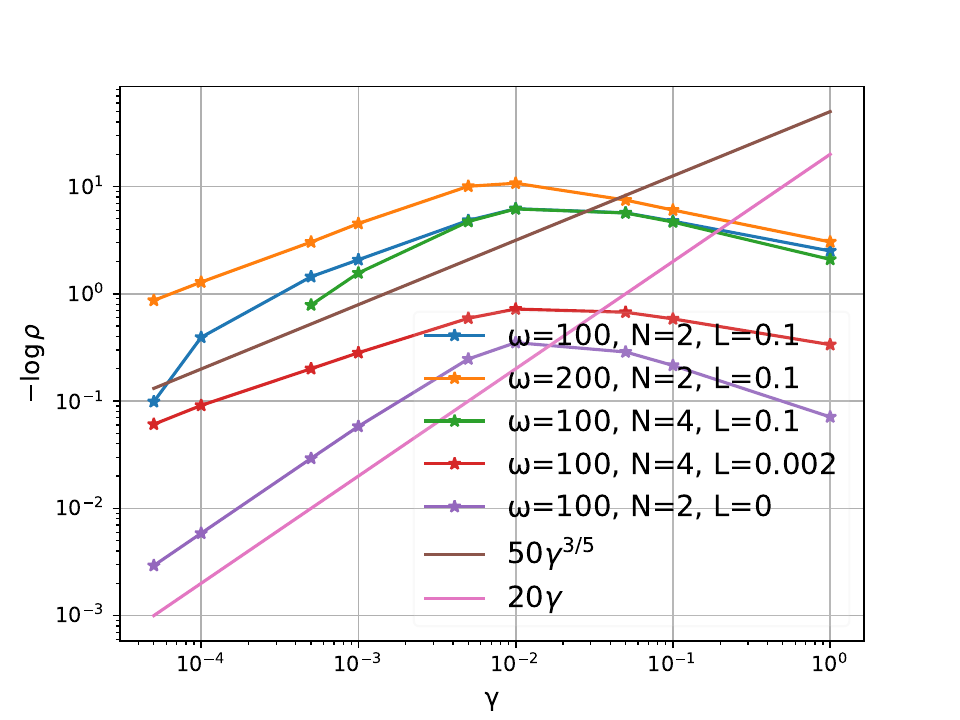}}
  \caption{Convergence factor dependence on $\gamma$ for the cavity
    with the operator $(1+ \I\omega \gamma)\Delta + \omega^2$}
  \label{fig13}
\end{figure}

\begin{figure}[t]
  \centering
  \mbox{\includegraphics[width=1.5in,trim=10 5 42 37,clip]{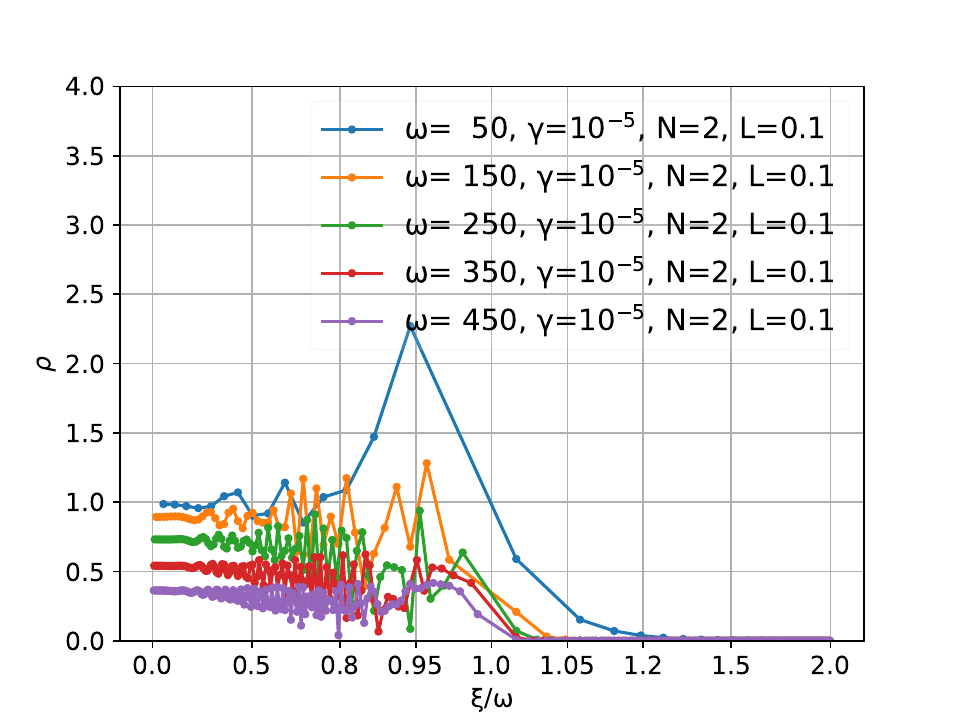}%
  \includegraphics[width=1.5in,trim=10 5 42 37,clip]{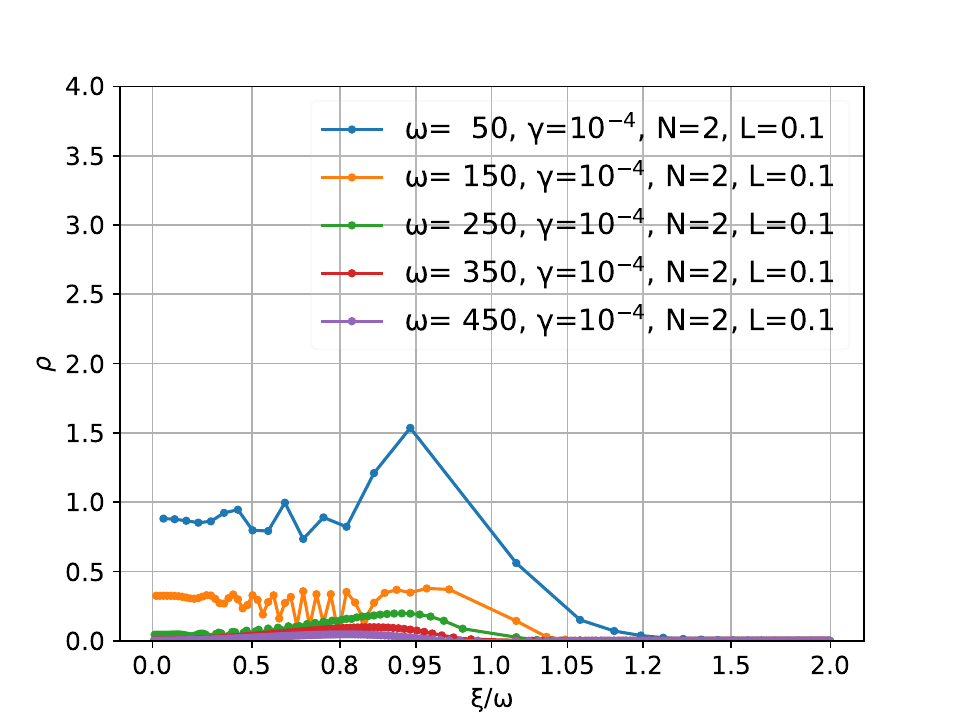}
  \includegraphics[width=1.5in,trim=10 5 45 37,clip]{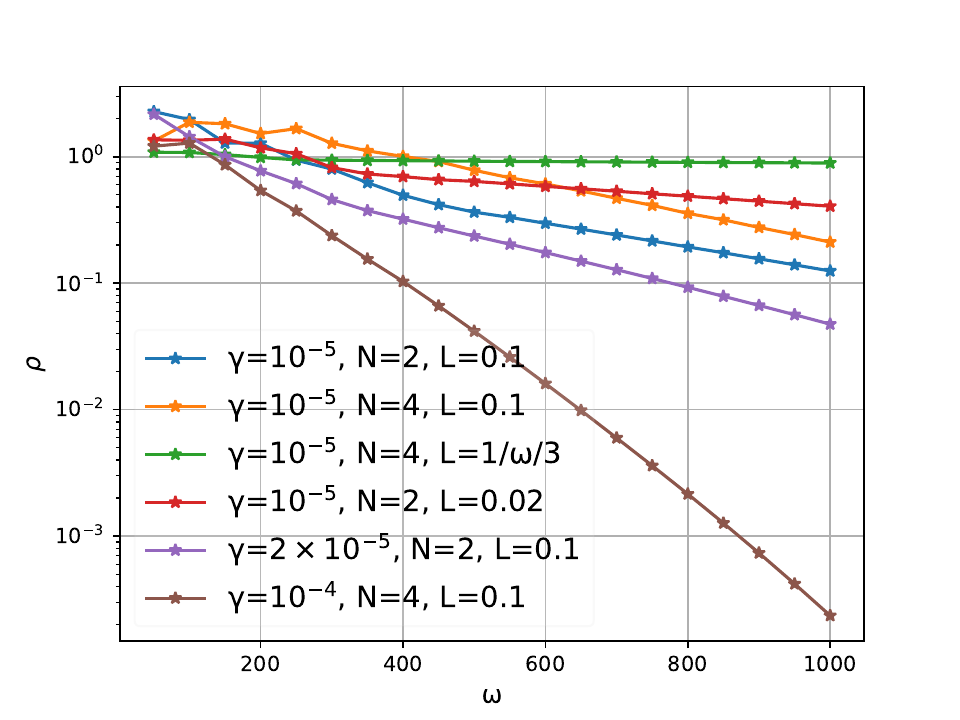}}

  \mbox{\includegraphics[width=1.5in,trim=10 5 42 37,clip]{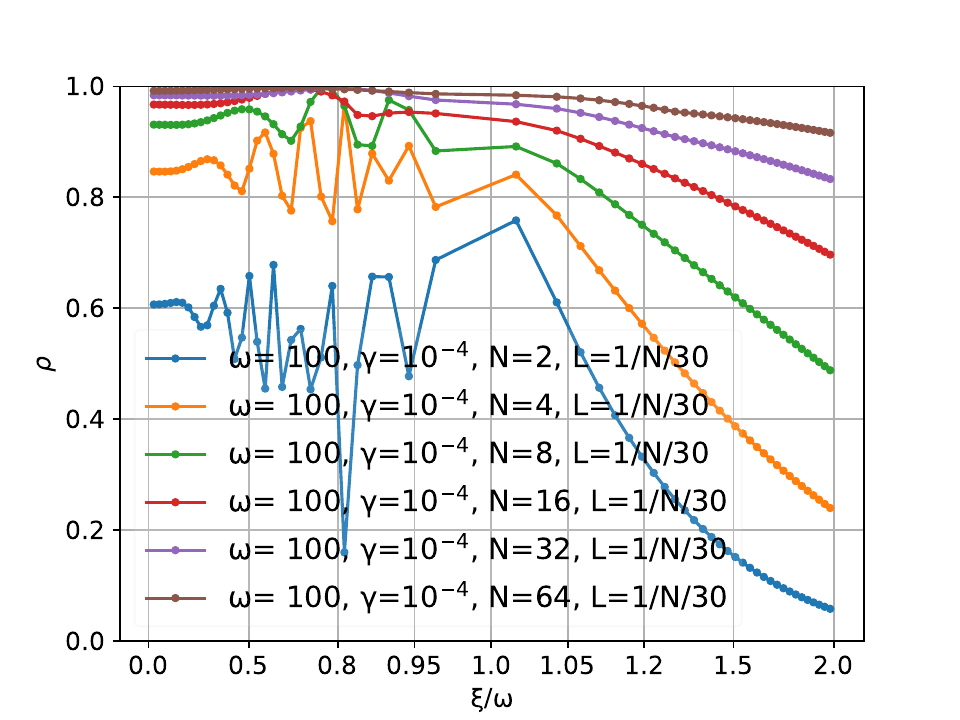}%
  \includegraphics[width=1.5in,trim=10 5 42 37,clip]{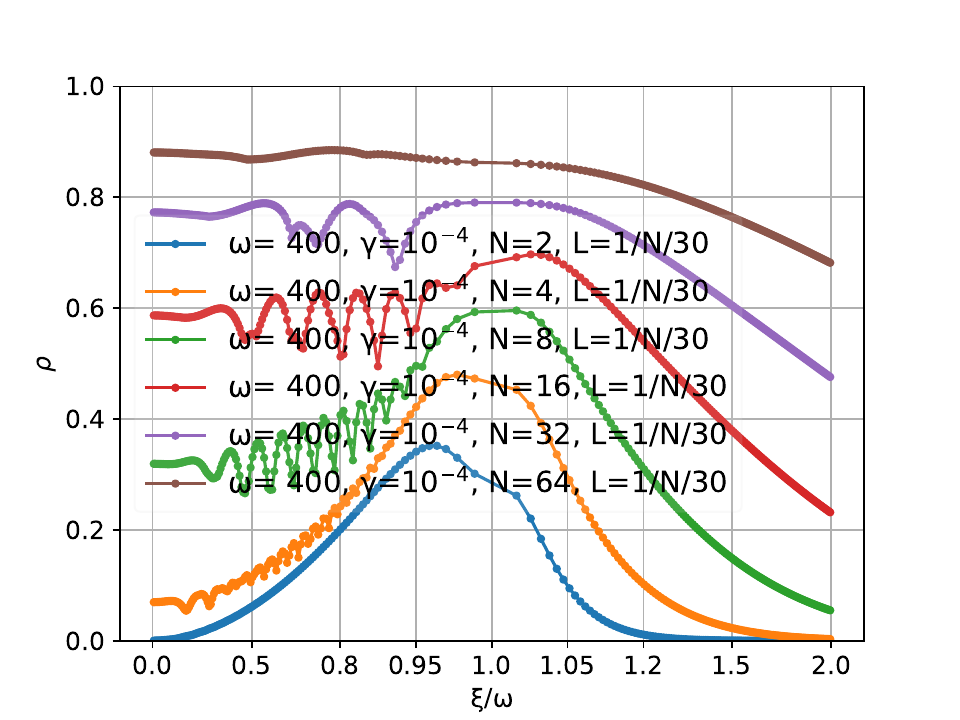}
  \includegraphics[width=1.5in,trim=10 5 45 37,clip]{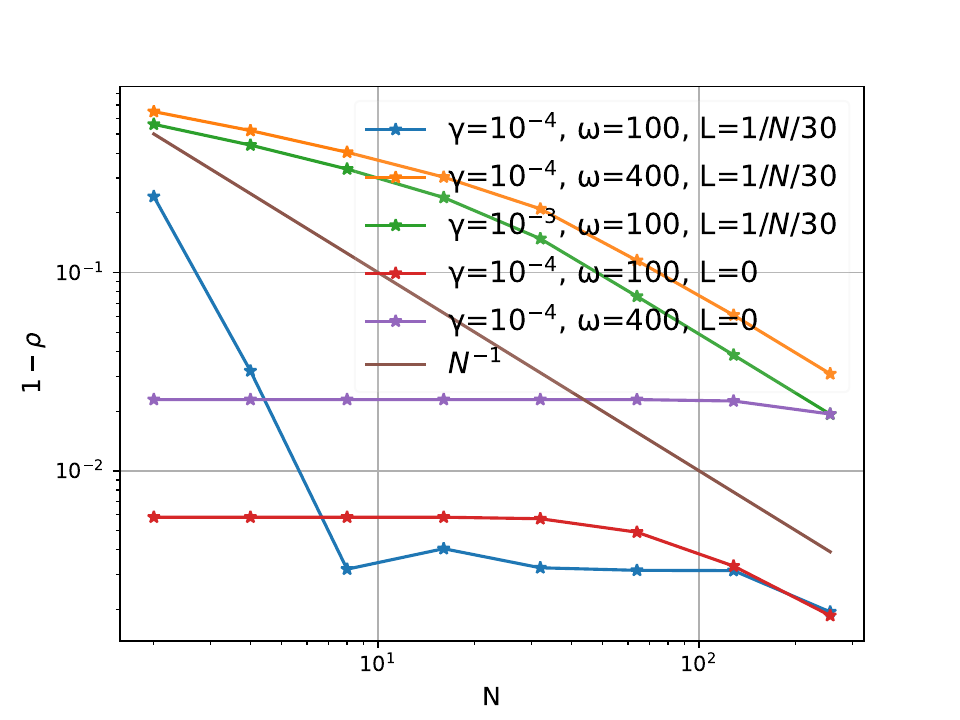}}

  \mbox{\includegraphics[width=1.5in,trim=10 5 42 37,clip]{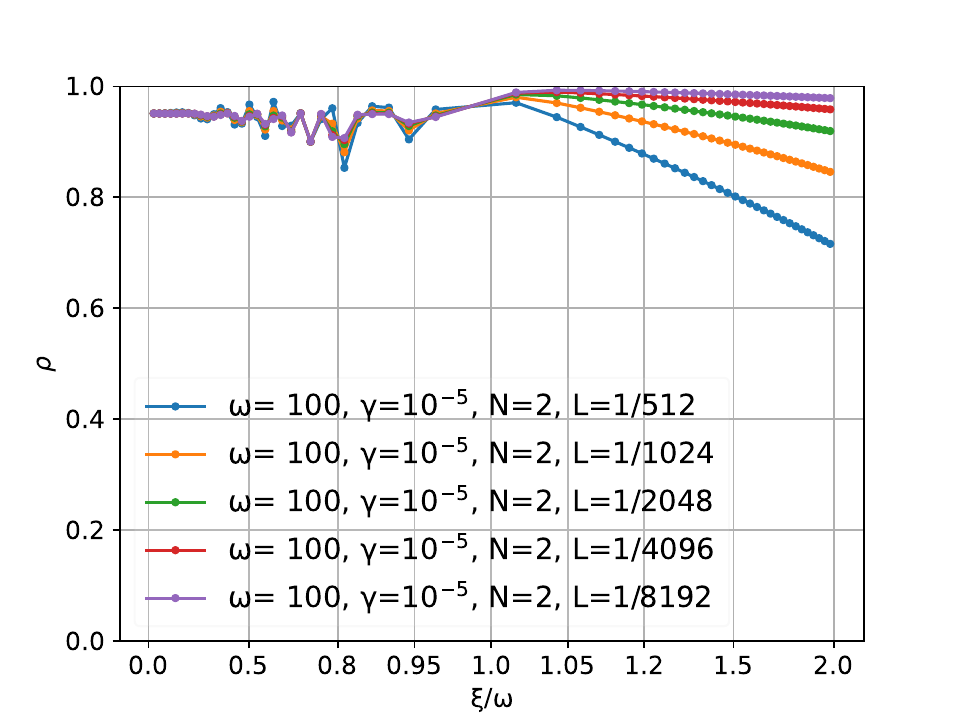}%
  \includegraphics[width=1.5in,trim=10 5 42 37,clip]{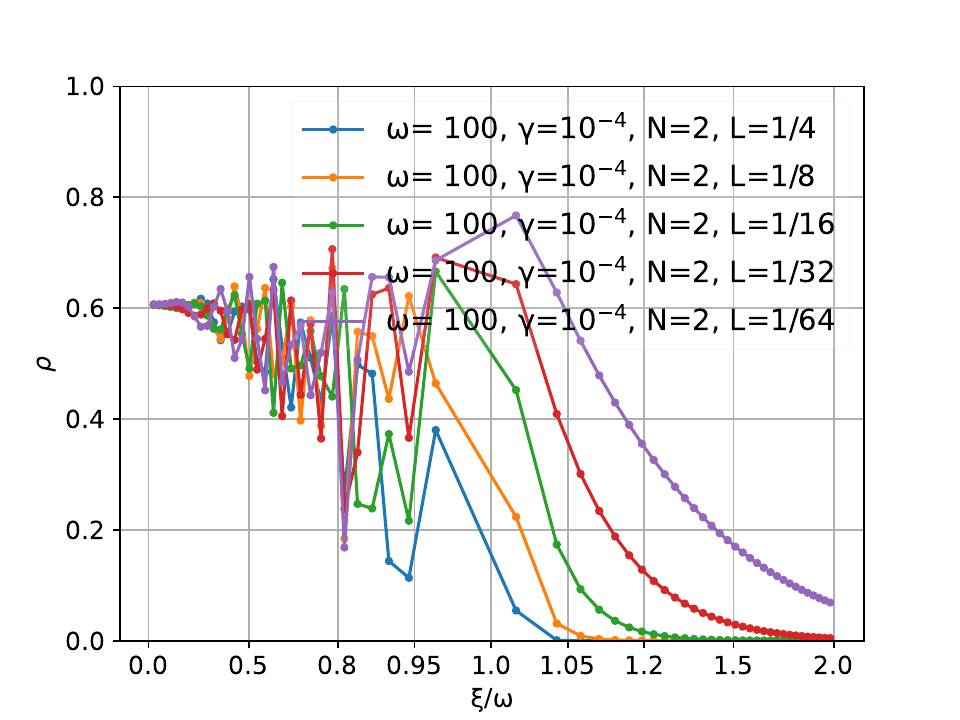}
  \includegraphics[width=1.5in,trim=10 5 45 37,clip]{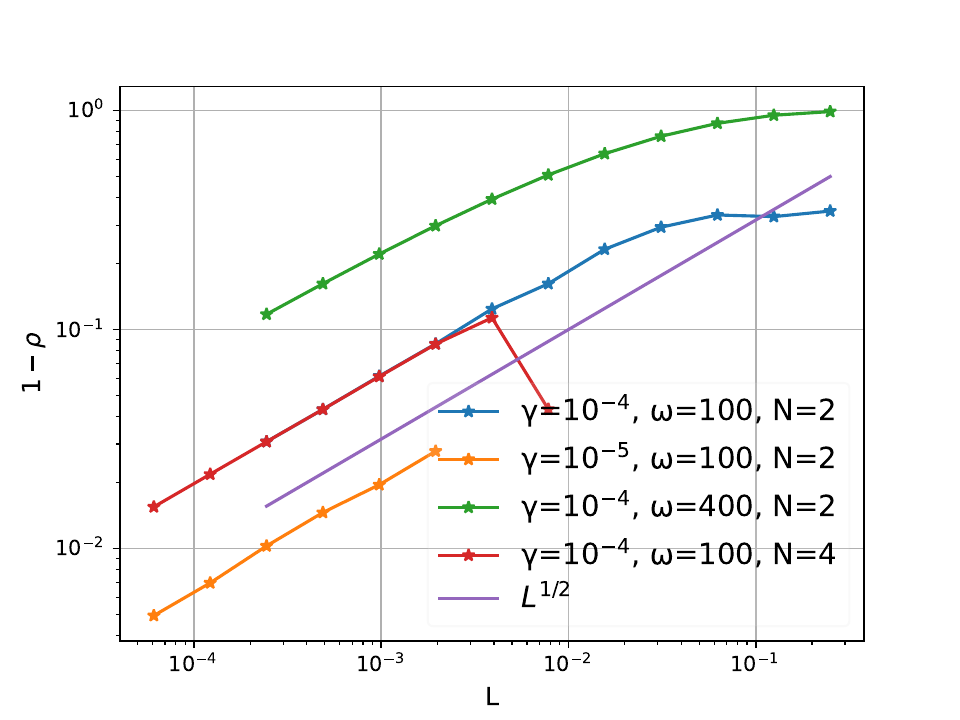}}
  \caption{Convergence factor dependence on $\omega$ (top), number of
    subdomains $N$ (middle) and overlap $L$ (bottom) for the cavity
    with the operator $(1+ \I\omega \gamma)\Delta + \omega^2$}
  \label{fig1416}
\end{figure}

We finally show the influence of the viscous damping coefficient $\gamma$ 
Fig.~\ref{fig13}.
We see that similarly to the first order damping, the behavior of the very hard closed cavity case
becomes comparable to the wave guide case, {especially at $\gamma=10^{-3}$, though,
  at smaller $\gamma$, the convergence at the space low frequency is slower than that in the wave
  guide.}  
The corresponding dependence on the
other parameters is shown in Fig.~\ref{fig1416}.

\bibliographystyle{spmpsci}
\bibliography{gander_zhang.bib}

\end{document}